\documentclass[sn-mathphys,Numbered]{sn-jnl}

\usepackage{graphicx}%
\usepackage{multirow}%
\usepackage{amsmath,amssymb,amsfonts}%
\usepackage{amsthm}%
\usepackage{mathrsfs}%
\usepackage[title]{appendix}%
\usepackage{xcolor}%
\usepackage{textcomp}%
\usepackage{manyfoot}%
\usepackage{booktabs}%
\usepackage{algorithm}%
\usepackage{algorithmicx}%
\usepackage{algpseudocode}%
\usepackage{listings}%

\theoremstyle{thmstyleone}%
\theoremstyle{thmstyletwo}%
\newtheorem{example}{Example}%

\theoremstyle{thmstylethree}%

\begin{document}

\title[A Numerically Robust and Stable Time-Space Pseudospectral Approach for Generalized Burgers-Fisher Equation]{A Numerically Robust and Stable Time-Space Pseudospectral Approach for Generalized Burgers-Fisher Equation}



\author[1]{\fnm{Harvindra} \sur{Singh}}\email{singh.harvindrahs@gmail.com}

\author*[1]{\fnm{L.K.} \sur{Balyan}}\email{lokendra.balyan@gmail.com}

\author[2]{\fnm{A.K.} \sur{Mittal}}\email{avinash.mittal10@gmail.com}

\author[1]{\fnm{Parul} \sur{Saini}}\email{iniaslurap1212@gmail.com}

\affil*[1]{\orgdiv{Discipline of Mathematics}, \orgname{PDPM IIITDM Jabalpur}, \orgaddress{\city{Jabalpur}, \postcode{482005}, \state{Madhya Pradesh}, \country{India}}}
\affil[2]{\orgdiv{Department of Mathematics}, \orgname{VIT Chennai}, \orgaddress{\city{Chennai}, \postcode{600127}, \state{Tamil Nadu}, \country{India}}}



\abstract{In this article, we present the time-space Chebyshev pseudospectral method (TS-CPsM) to approximate a solution to the generalised Burgers-Fisher (gBF) equation. The Chebyshev-Gauss-Lobatto (CGL) points serve as the foundation for the recommended method, which makes use of collocations in both the time and space directions.
Further, using a mapping, the non-homogeneous initial-boundary value problem is transformed into a homogeneous problem and a system of algebraic equations is obtained.
The numerical approach known as Newton-Raphson is implemented in order to get the desired results for the system. The proposed method's stability analysis has been performed. Different researchers' considerations on test problems have been explored to illustrate the robustness and practicality of the approach presented.
The approximate solutions we found using the proposed method are highly accurate and are significantly better than the existing results.  
}




\maketitle

\section{Introduction}\label{sec1}

The nonlinearity of equations has been a puzzling problem for researchers for decades.
Numerous studies have been conducted and are being conducted to address this concern and help science move forward.
Some of the most popular non-linear equations are Burgers and Fishers equations which are of great importance due to their wide range of applications in mathematical physics. 
In 1915, Harry Bateman initially proposed the Burgers equation 
\begin{align}\label{Eq.1.1}
U_{\tilde{t}} + U U_{\tilde{\eta}} - \mu U_{\tilde{\eta} \tilde{\eta}}=0
\end{align}
 which was further investigated in 1948 by Johannes Martinus Burgers~\cite{bateman1915some},~\cite{BURGERS1948171}. The applications of Burgers equation can be seen in the context of traffic flow~\cite{greenshields1935study},
 fluid mechanics~\cite{murray1973burgers}, gas dynamics~\cite{hirsh1975higher}, and quantum field~\cite{yepez2002efficient}.
  Fishers equation
 \begin{align}\label{Eq.1.2}
 U_{\tilde{t}} -\mu   U_{\tilde{\eta} \tilde{\eta}}=\sigma_2 U(1-U)
 \end{align}
  referred to in certain circles as KPP(Kolmogorov-Petrovsky-Piskunov) equation or Fisher-KPP equation is a kind of equation that takes into account both reactions and diffusions which was introduced by R. A. Fisher in 1937~\cite{fisher1937wave}.
 The Fishers equation (reaction-diffusion equation) appears in the modelling of a wide range of scientific and physical phenomena, such as  Brownian motion~\cite{bramson1978maximal}, ecology~\cite{britton1986reaction}, combustion~\cite{frank2015diffusion}, nuclear reactor theory~\cite{canosa1969diffusion}, population dynamics~\cite{fisher1937wave}, physiology, etc.
A nonlinear equation is obtained by combining these two equations \ref{Eq.1.1} and \ref{Eq.1.2} called Burgers-Fisher equation~\cite{ablowitz1987topics} and is given by
\begin{align} \label{Eq.3}
U_{\tilde{t}} + \sigma_1 UU_{\tilde{\eta}} - \mu U_{\tilde{\eta} \tilde{\eta}} = \sigma_2 U(1-U). 
\end{align} 
This equation consists of the effects of non-linear advection, linear diffusion, and non-linear logistic reactions.  This equation is a preliminary model that may be used to describe the interaction between the reaction mechanism, the diffusive transport, and the convective effect~\cite{wang1990exact}.
The Burgers-Fisher equation is useful for solving a wide variety of issues involving applied mathematics and physics. Its aplications include gas dynamics, traffic flow, number theory, elasticity, heat conduction, and financial mathematics~\cite{chandrasekaran1996painleve,chen2004new,lu2007some}. 
The generalized Burgers-Fisher equation, abbreviated as the gBF equation, is the most comprehensive representation of the Burgers-Fisher equation and is given by  
\begin{align} \label{Eq.1.4}
U_{\tilde{t}} + \sigma_1 U^{\delta} U_{\tilde{\eta}} -\mu U_{\tilde{\eta} \tilde{\eta}} = \sigma_2 U(1-U^{\delta}), ~~~~0 \leq \tilde{\eta} \leq 1,~~~~ \tilde{t} \geq 0,
\end{align} 
where $\sigma_1 \geq 0$,  $\mu \geq 0$, and $\sigma_2 \geq 0$ represent the advection coefficient of the non-linear advection term, the diffusion coefficient of linear diffusion term, and the reaction coefficient of non-linear reaction term, respectively. 
Literature shows a very high range of applications of Eq.~\ref{Eq.1.4} in the field of physics, including chemical physics, plasma physics, and fluid physics. Applications are not limited and are reported in various other fields too, such as nonlinear optics, capillary gravity waves, etc.\newline 
Due to its extreme non-linearity, the generalized Burgers-Fisher equation has served as a model for numerous academics to try out new numerical and analytical techniques.
In 2004, an explicit solution was given by Kaya and El-Sayed via a numerical simulation for the gBF equation~\cite{kaya2004numerical}. Further, a b-spline quasi-interpolation method was proposed by Zhu and Kang~\cite{zhu2010numerical} in 2010 to approximate the gBF equation. 
Cubic b-spline collocation and fourth-order b -spline collocation method~\cite{mittal2015numerical},~\cite{singh2020fourth}, homotopy perturbation method(HPM)~\cite{rashidi2009explicit}, Adomian decomposition method(ADM)~\cite{ismail2004adomian}, collocation-based radial basis functions method (CBRBF)~\cite{khattak2009computational}, Exp-function method~\cite{xu2010application}, Lattice Boltzmann model~\cite{zhang2010lattice}, etc. have been used by several researchers to find the numerical as well as the analytical solutions of the gBF equation.\newline 
Further, a compact finite difference method and a polynomial-based quadrature method were proposed by Sari $et$ $al.$ in order to approximate the gBF equation~\cite{sari2010compact},~\cite{sari2011differential}.
A local discontinuous Galerkin (LDG) methods and optimal homotopy asymptotic method (OHAM) have been employed by Zhang $et$ $al.$~\cite{zhang2012local} and Nawaz $et$ $al.$~\cite{nawaz2013application}, respectively. Galerkin finite element method was employed by Yadav and Jiwari~\cite{yadav2017finite} to approximate and analyze the Burgers-Fisher equation.
Malik $et$ $al.$~\cite{malik2015numerical} numerically simulate the gBF equation using  Exp-function method.\newline
The most recent work on the gBF equation is referred to~\cite{nadeem2019modified,loyinmi2020exact,akinfe2021solitary,mohanty2020high,kumar2019stability,li2019geometric,onyejekwe2020numerical,zhang2021global}.
 The research articles~\cite{mittal2021numerical, mittal2020spectrally,balyan2020stability,saini2022comparative,saini2022modification} and the references therein provide further information on the pseudospectral approach for both smooth and non-smooth functions.\newline
Within the past few years, the spectral method has gained popularity as a reliable scientific computing tool for numerically solving linear and strongly non-linear partial differential equations~\cite{hesthaven2007spectral,boyd2001chebyshev,trefethen2000spectral,canuto2012others}. The method's high-order accuracy has led to its widespread use. Interpolation and non-interpolation are the broad categories into which the spectral method has been classified. The pseudospectral method, also known as the spectral collocation method, depends on the interpolation methodology.
The pseudospectral method forms a class of techniques that are efficient and highly accurate and having the potential for rapidly convergent approximations and low error estimation. This article will examine the gBF equation.
In 2006, an approximation of the gBF equation using a modified version of the pseudospectral approach was given by Javidi~\cite{javidi2006modified} with tremendous accuracy in which the pseudospectral method was applied in the space direction only.
For the numerical approximation of the gBF equation, we offer a spectrally accurate enhanced version of the Chebyshev pseudospectral approach that takes place in time as well as in space, and the results obtained using the suggested approach have proved that time is just as significant as it utilizes in space.
 To the best of the author's knowledge, this is the first work to employ the pseudospectral technique for the approximation of the gBF equation in both time and space.\newline 
The rest of the article is arranged as follows: We provide an in-depth explanation of the suggested approach as well as a clear image of how the discretization process works in Sect.~\ref{sec2}. Analysis of the method's stability for the gBF equation is conducted in Sect.~\ref{sec3}.  The accuracy and efficiency of the suggested technique are shown in Sect.~\ref{sec4}. Lastly, the observations and the concluding remarks are given in Sect.~\ref{sec5}.

    \section{Chebyshev pseudospectral discretization of the gBF equation in both time and space}\label{sec2}
    In this section, we discuss the Chebyshev pseudospectral approach to approximate the gBF equation in both time and space.
We look for a pseudospectral approximation $U(\eta, t)$ in the form of a finite linear combination of a product of Chebyshev polynomials with spectral coefficients such that
    \begin{equation}\label{3.1}
    I_N U= \sum_{\alpha=0}^{N}\sum_{\beta=0}^{N} \phi_{\alpha}(\eta)  \phi_{\beta}(t) \Theta_{\alpha\beta},
    \end{equation}
    where, in space and time directions, Chebyshev polynomials are $\phi_{\alpha}(\eta)$ and $\phi_{\beta}(t)$ respectively and spectral coefficients $\Theta_{\alpha\beta}$ are given by 
     \begin{equation}
      \Theta_{\alpha\beta}= \frac{1}{\hat{h}_{\alpha}  \hat{h}_{\beta}} \sum_{l=0}^{N}\sum_{m=0}^{N} \phi_{\alpha}(\eta_{l}) \phi_{\beta}(t_{m}) w_{l} w_{m} U(\eta_{\alpha},t_{\beta}),      
     \end{equation}
     where $\hat{h}_{p} = \frac{\epsilon_p}{2} \pi$ is a discrete normalizing factor, $\epsilon_0 = 2$ and $\epsilon_p = 1$ for $p\geq 1$.
   Where $p = \{ \alpha, \beta \}$ and $w$ denotes the weight function.\newline
     Also, we can write 
\begin{align}\label{Eq.2.3}
\nonumber I_N U &= \sum_{\alpha=0}^{N}\sum_{\beta=0}^{N} \phi_{\alpha}(\eta) \phi_{\beta}(t) \Theta_{\alpha\beta}, \\
 & = \nonumber \sum_{\alpha=0}^{N}\sum_{\beta=0}^{N} \phi_{\alpha}(\eta)  \phi_{\beta}(t) \frac{1}{\hat{h}_{\alpha}\hat{h}_{\beta}} \sum_{l=0}^{N} \sum_{m=0}^{N} \phi_{\alpha}(\eta_{l})  \phi_{\beta}(t_{m}) w_{l} w_{m} U(\eta_{\alpha},t_{\beta}),\\
 & \nonumber = \nonumber \sum_{\alpha=0}^{N}\sum_{\beta=0}^{N}  \left(\frac{1}{\hat{h}_{\alpha}\hat{h}_{\beta}} \sum_{l=0}^{N} \sum_{m=0}^{N}    w_{l}   \phi_{\alpha}(\eta_{l}) \phi_{\alpha}(\eta)   w_{m} \phi_{\beta}(t_{m}) \phi_{\beta}(t) \right) U (\eta_{\alpha},t_{\beta}),\\
 &  = \sum_{\alpha=0}^{N} \sum_{\beta=0}^{N} \mathcal{L}_{\alpha}(\eta)  \mathcal{L}_{\beta}(t) U(\eta_{\alpha},t_{\beta}), 
 \end{align}
where $\mathcal{L}$ denotes the Lagrange polynomials which can be written as 
\begin{align*}
 & \mathcal{L}_p(z) = \frac{1}{\hat{h}_p} \sum_{r=0}^{N} w_{r} \phi_{p} (z_{r}) \phi_{p} (z) = \frac{(-1)^{p+1} (1-z^2)}{c_{p} N^2(z-z_{p})} \frac{\partial}{\partial z} \phi_N(z),
\end{align*}
where, $p = \{\alpha, \beta \}$,$r =\{l, m\}$, $z = \{\eta,t\}$,
and
$$c_{p} = \begin{cases} 2, ~~~~~~~~ p = 0 ~~ \text{and}~~ N,\\ 1, ~~~~~~~~ 0 < p <N.
    \end{cases}$$   
Here, the Lagrange polynomial satisfies the delta function i.e. $\mathcal{L}_{p}(z_{\beta}) = \delta_{p \beta}$. The differential matrix $\mathcal{D}$ of first order is given by,
\begin{equation} 
\mathcal{D}_{\alpha \beta}=\begin{cases} -\frac{2N^2 + 1}{6},~~~~ \alpha=\beta=0,\\ \frac{c_\alpha}{c_\beta} \frac{(-1)^{\alpha + \beta}}{z_\alpha - z_\beta}, ~~~~ \alpha \neq \beta,\\
-\frac{z_\alpha}{2(1-z_\alpha^2)},~~~~\alpha=\beta \in [1,...,N-1],\\
\frac{2N^2 + 1}{6},~~~~ \alpha=\beta=N
\end{cases}
\end{equation}
It is possible to write the spectral approximation given in Eq.~\ref{Eq.2.3} as a direct product
\begin{align}\label{Eq.2.5}
\nonumber I_N U & = [\mathcal{L}_0(\eta) \mathcal{L}_0(t),...\mathcal{L}_0(\eta)\mathcal{L}_N(t),...,\mathcal{L}_{N}(\eta)\mathcal{L}_0(t),...,\mathcal{L}_{N}(\eta) \mathcal{L}_N(t)] \boldsymbol{U} ,\\
 & = \nonumber
  [\left\lbrace \mathcal{L}_0(\eta),...,\mathcal{L}_N(\eta) \right\rbrace \otimes \left\lbrace \mathcal{L}_{0}(t),..., \mathcal{L}_N(t)\right\rbrace] \boldsymbol{U},  \\ 
 &   =  \nonumber
    (\psi_{[0:N]}(\eta) \otimes \psi_{[0:N]}(t))^T \boldsymbol{U},\\
 &   =  
 (K_{[0:L]} (\eta,t))^T \boldsymbol{U}.
\end{align}
 where 
  \begin{equation}
\nonumber \boldsymbol{U} = [U_{00},...,U_{0N}|U_{10},...,U_{1N}|...|U_{N0},...,U_{NN}]^T,
\end{equation}
\begin{equation}
\nonumber \psi_{[0:N]}(z) = [\mathcal{L}_0(z),...,\mathcal{L}_N(z)]^T,~~~~  \text{and}~~~~ 
 L = (N+1) \times (N+1) - 1.
 \end{equation} 
    The first and second-order derivatives of $U$ with respect to '$\eta$' are now determined using Eq.~\ref{Eq.2.5}
  \begin{equation}
  \frac{\partial U (\eta,t)}{\partial \eta} = (K_{[0:L]} (\eta,t))^T (\mathcal{D}^{'}_{[0:N,0:N]} \otimes I_{N+1}) \boldsymbol{U},
  \end{equation}  
    \begin{equation}
  \frac{\partial^2 U (\eta,t)}{\partial \eta^2} = (K_{[0:L]} (\eta,t))^T (\mathcal{D}^{''}_{[0:N,0:N]} \otimes I_{N+1}) \boldsymbol{U},
  \end{equation}  
  where $\mathcal{D}^{'}_{[0:N,0:N]}$ and $\mathcal{D}^{''}_{[0:N,0:N]}$ denotes the first and the second order differential matrices. 
  Moreover, the first derivative with respect to '$t$' is given by 
    \begin{equation}
    \frac{\partial U(\eta,t)}{\partial t} = (K_{[0:L]} (\eta,t))^T (I_{N+1} \otimes \mathcal{D}_{[0:N,0:N]}^t) \boldsymbol{U}.
    \end{equation} 
  Now, we will discretize the gBF equation based on discussed pseudospectral method
and in order to apply the method the idea is to convert the space and the time intervals into $[-1,1]$ and therefore to serve the purpose we use the following transformation for any $\tilde{v} \in [I_L, I_R] $
$$\tilde{v} \rightarrow \frac{I_R-I_L}{2}v + \frac{I_R+I_L}{2}.$$
     By using the transformations in both space and time, we obtain the following generalized Burgers-Fisher equation in new time and space intervals  
       \begin{align}\label{Eq.2.9}
      U_t + A (U^{\delta+1})_{\eta} - B U_{\eta \eta} - C (U-U^{\delta + 1 }) = 0,
       \end{align}
       with the initial condition
       \begin{align}\label{Eq.2.10}
       U(\eta,-1) = h (\eta),~~~~~\eta\in[-1,1],
       \end{align} 
       and the associated boundary conditions
         \begin{align}\label{Eq.2.11}
         U(-1,t) = g _1(t),~~  U(1,t) = g_2(t),  ~~~t\in[-1,1],
         \end{align} 
 where
 
 \begin{align}\label{Eq.2.12}
  A=\frac{\sigma_1 T}{(\delta  + 1)(\eta_R - \eta_L)}, ~B=\frac{2T}{(\eta_R - \eta_L)^2}, ~\text{and}~  C=\frac{\sigma_2 T}{2}.
 \end{align}
              Further, to make it easier to solve, a mapping is presented for changing a system's initial and boundary values from non-homogeneous to homogeneous,
          
          \begin{align}\label{Eq.2.13}
  \nonumber  \Omega(\eta,t) &= \frac{(1-t)}{2} h(\eta) + \frac{(1-\eta)}{2} g_1(t) + \frac{(1+\eta)}{2} g_2(t) - \frac{(1-t)(1-\eta)}{4}g_1(-1)\\
   & - \frac{(1-t)(1+\eta)}{4}g_2(-1),
          \end{align}
        where at the corners of the domain the initial-boundary conditions are supposed to satisfy
         $h (-1) = g_1 (-1)$
          and $h (1) = g_2 (-1).$
          Now, we have transformed our non-homogeneous initial-boundary value problem into homogeneous one. Let us define a variable $V$, such that
          \begin{align}\label{Eq.2.14}
          U(\eta, t) = V(\eta, t) + \Omega(\eta, t).
          \end{align}
          Using Eq.~\ref{Eq.2.14}, Eqs.~\ref{Eq.2.9}-\ref{Eq.2.12} can be written as
          \begin{align}\label{Eq.2.15}
       \nonumber &  (V(\eta, t) + \Omega(\eta, t))_t + A [(V(\eta, t) + \Omega(\eta, t))^{\delta+1}]_\eta - B (V(\eta, t) + \Omega(\eta, t))_{\eta \eta}\\
       & - C[(V(\eta, t)
               + \Omega(\eta, t))-(V(\eta, t) + \Omega(\eta, t))^{\delta + 1}] = 0,
           \end{align}
           with the following homogeneous initial and boundary conditions
           \begin{align}
           &  V(\eta, -1) = 0,~~~~~\eta\in[-1,1],\\
          &\label{Eq.2.17} V(-1,t) = 0,~~  V(1,t) = 0,  ~~~t\in[-1,1],
           \end{align}   
      Using Binomial expansion in Eq.~\ref{Eq.2.15}, we get 
          \begin{align}\label{Eq.2.18}
       \nonumber &  (V(\eta, t) + \Omega(\eta, t))_t + A [~^{\delta+1}c_0 V^{\delta+1}(\eta,t) \Omega^{0}(\eta,t)~+~^{\delta+1}c_1 V^{\delta}(\eta,t) \Omega^{1}(\eta,t) +...\\
       & \nonumber +~^{\delta+1}c_{\delta} V^{1}(\eta,t) \Omega^{\delta}(\eta,t) +~ ^{\delta+1}c_{\delta+1} V^{0}(\eta,t) \Omega^{\delta+1}(\eta,t)]_\eta - B (V(\eta, t) + \Omega(\eta, t))_{\eta \eta}\\
      & \nonumber - C(V(\eta, t) + \Omega(\eta, t)) + C[~^{\delta+1}c_0 V^{\delta+1}(\eta,t) \Omega^{0}(\eta,t) ~+~^{\delta+1}c_1 V^{\delta}(\eta,t) \Omega^{1}(\eta,t) +...\\
      & +~^{\delta+1}c_{\delta} V^{1}(\eta,t) \Omega^{\delta}(\eta,t)+~ ^{\delta+1}c_{\delta+1} V^{0}(\eta,t) \Omega^{\delta+1}(\eta,t)] = 0.
           \end{align}
     Now, as described in detail at the beginning of this section, we look for an approximation that can be defined in vector-vector product form, therefore, the spectral approximation $\boldsymbol{V}(\eta,t)$ is defined as   
                  \begin{align}
                    \nonumber & \boldsymbol{V}(\eta,t)=(\psi_{[1:N-1]}(\eta) \otimes \psi_{[1:N]}(t))^T V, \\
                 & \label{Eq.2.19} \boldsymbol{V}(\eta,t) = (K_{[1:(N-1)N]} (\eta,t))^T V, 
                  \end{align}
                  where $V$ is a $(N-1) \times N $ vector, given by
                  \begin{align*}
                  V = [V_{11},...,V_{1N}|V_{21},...,V_{2N}|...|V_{(N-1)1},...,V_{(N-1)N}]^T.
                  \end{align*}
                  Additionally, we can define the following derivatives of $\boldsymbol{V}(\eta,t)$ 
                  \begin{align}\label{Eq.2.20}
                     \boldsymbol{V}_{t}(\eta,t) = (K_{[1:(N-1)N]} (\eta,t))^T ( I_{N} \otimes  \mathcal{D}^t ) V,
                  \end{align} 
                  \begin{align}\label{Eq.2.21}
                      \boldsymbol{V}_{\eta}(\eta,t) = (K_{[1:(N-1)N]} (\eta,t))^T (\mathcal{D}^{\eta'}  \otimes I_N) V,
                  \end{align}
                    \begin{align}\label{Eq.2.22}
                      \boldsymbol{V}_{\eta \eta}(\eta,t) =  (K_{[1:(N-1)N]} (\eta,t))^T (\mathcal{D}^{\eta''}  \otimes I_N) V,
                  \end{align}
   where $\mathcal{D}^{t}=\mathcal{D}^{t}_{[1:N,1:N]}$ is the first order differential matrix w.r.t. '$t$', and $\mathcal{D}^{\eta'}=\mathcal{D}^{\eta'}_{[1:N-1,1:N-1]}$ and $\mathcal{D}^{\eta''}=\mathcal{D}^{\eta''}_{[1:N-1,1:N-1]}$  are the first and the second order differential matrices with respect to '$\eta$', respectively.\newline
      By using Eqs.~\ref{Eq.2.19}-\ref{Eq.2.22} in Eq.~\ref{Eq.2.18}, we obtain
     \begin{align}\label{Eq.2.23}
     \nonumber & (K_{[1:(N-1)N]} (\eta,t))^T \Bigl\{ (I_{N-1} \otimes \mathcal{D}^t)V+
      A(\mathcal{D}^{\eta'}  \otimes I_N) V^{\delta +1} + A~^{\delta+1}c_1~  [(\text{diag}(\Omega_{\eta}(\eta,t))) \\
      & \nonumber  + \text{diag}(\Omega(\eta,t))(\mathcal{D}^{\eta} \otimes I_N)] V^{\delta}  +...+ A~^{\delta+1}c_{\delta}~  [(\text{diag}((\Omega^{\delta}(\eta,t))_{\eta}))\\
            & \nonumber + \text{diag}(\Omega^{\delta}(\eta,t))(\mathcal{D}^{\eta} \otimes I_N)] V  - B (  \mathcal{D}^{\eta''} \otimes I_N)V   + C[-(I_{N-1} \otimes I_N)V\\
                  & \nonumber + (I_{N-1} \otimes I_N)V^{\delta+1} +  ~^{\delta+1}c_1~  (\text{diag}(\Omega(\eta,t))) V^{\delta} +...+  ~^{\delta+1}c_{\delta}~  (\text{diag}(\Omega^{\delta}(\eta,t))) V]   \Bigr\}\\
                  &  + \Omega_t(\eta,t) + A(\Omega^{\delta+1} (\eta,t))_\eta - B \Omega_{\eta \eta}(\eta,t) -C\Omega(\eta,t)  + C \Omega^{\delta+1}(\eta,t)= 0.
                \end{align}           
                 Further, in Eq.~\ref{Eq.2.23}, using Chebyshev Gauss Lobbatto points in $\eta$ and $t$ directions such that $(\eta,t) = (\eta_i,t_j),~~i=1,..., N-1, ~~\text{and}~~j=1,..., N$, and using the fact that the Lagrange polynomial satisfies the properties of a delta function and hence $(K_{[1:(N-1)N]} (\eta_i,t_j))^T =  [\psi_{[1:N-1]}(\eta_i) \otimes \psi_{[1:N]}(t_j)]^T$ gives a vector $V_r^T$ of order $(N-1) \times N$ in which except the $r^{th}$ element rest are zero.
                 Therefore, we obtain 
     \begin{align}\label{Eq.2.24}
     \nonumber &   (I_{N-1} \otimes \mathcal{D}^t)V+
      A(\mathcal{D}^{\eta'}  \otimes I_N) V^{\delta +1} + A~^{\delta+1}c_1~  [(\text{diag}(\Omega_{\eta}(\eta_i,t_j)))\\
            & \nonumber + \text{diag}(\Omega(\eta_i,t_j))(\mathcal{D}^{\eta} \otimes I_N)] V^{\delta}  +...+ A~^{\delta+1}c_{\delta}~  [(\text{diag}(\delta \Omega^{\delta-1}(\eta_i,t_j) \Omega_{\eta}(\eta_i,t_j)))\\
                  & \nonumber + \text{diag}(\Omega^{\delta}(\eta_i,t_j))(\mathcal{D}^{\eta} \otimes I_N)] V - B (  \mathcal{D}^{\eta''} \otimes I_N)V  \\
       & \nonumber + C[-(I_{N-1} \otimes I_N)V + (I_{N-1} \otimes I_N)V^{\delta+1} +  ~^{\delta+1}c_1~  (\text{diag}(\Omega(\eta_i,t_j))) V^{\delta}\\
             & \nonumber +...+  ~^{\delta+1}c_{\delta}~  (\text{diag}(\Omega^{\delta}(\eta_i,t_j))) V ]  + \Omega_t(\eta_i,t_j) + (\delta+1)A(\Omega^{\delta} (\eta_i,t_j)) (\Omega_\eta(\eta_i,t_j))\\
             & - B \Omega_{\eta \eta}(\eta_i,t_j)  -C\Omega(\eta_i,t_j)  + C \Omega^{\delta+1}(\eta_i,t_j)= 0.
                \end{align} 
                In conclusion, the equation system~\ref{Eq.2.24} can be converted into the matrix form as follows  
     \begin{align}\label{Eq.2.25}
     \nonumber &   (I_{N-1} \otimes \mathcal{D}^t)V+
      A(\mathcal{D}^{\eta'}  \otimes I_N) V^{\delta +1} + A~^{\delta+1}c_1~  [(\text{diag}(\textbf{Q})) + \text{diag}(\textbf{P})(\mathcal{D}^{\eta} \otimes I_N)] V^{\delta} \\
      & \nonumber +...+ A~^{\delta+1}c_{\delta}~  [(\text{diag}(\delta \textbf{P}^{\delta-1} \textbf{Q}) + \text{diag}(\textbf{P}^{\delta})(\mathcal{D}^{\eta} \otimes I_N)] V - B (  \mathcal{D}^{\eta''} \otimes I_N)V  \\
       & \nonumber + C [-(I_{N-1} \otimes I_N)V + (I_{N-1} \otimes I_N)V^{\delta+1} +  ~^{\delta+1}c_1~  (\text{diag}(\textbf{P})) V^{\delta} \\
       & +...+  ~^{\delta+1}c_{\delta}~  (\text{diag}(\textbf{P}^{\delta})) V ] + \textbf{S} + (\delta+1)A\textbf{P}^{\delta} \textbf{Q} - B \textbf{R}  -C\textbf{P}  + C \textbf{P}^{\delta+1}= 0,
                \end{align}               
                 where
                \begin{align}\label{Eq.2.26}
                \nonumber &V = [V_{11},...,V_{1N}|V_{21},...,V_{2N}|...|V_{(N-1)1},...,V_{(N-1)N}]^T,\\
                & \nonumber \textbf{P} = [\Omega(\eta_1,t_1),...,\Omega(\eta_1,t_N)|\Omega(\eta_2,t_1),...,\Omega(\eta_2,t_N)|... |\Omega(\eta_{N-1},t_1), ...,\Omega(\eta_{N-1},t_N)]^T,\\
                & \nonumber \textbf{Q} = [\Omega_{\eta}(\eta_1,t_1),...,{\Omega_{\eta}(\eta_1,t_N)}|\Omega_{\eta}(\eta_2,t_1),...,\Omega_{\eta}(\eta_2,t_N)|... |\Omega_{\eta}(\eta_{N-1},t_1),...,\Omega_{\eta}(\eta_{N-1},t_N)]^T,\\
                & \nonumber \textbf{R} = [\Omega_{\eta \eta}(\eta_1,t_1),...,\Omega_{\eta \eta}(\eta_1,t_N)|\Omega_{\eta \eta}(\eta_2,t_1),...,\Omega_{\eta \eta}(\eta_2,t_N)|... |\Omega_{\eta \eta}(\eta_{N-1},t_1),...,\Omega_{\eta \eta}(\eta_{N-1},t_N)]^T,\\
                & \textbf{S} = [\Omega_{t}(\eta_1,t_1),...,{\Omega_{t}(\eta_1,t_N)}|\Omega_{t}(\eta_2,t_1),...,\Omega_{t}(\eta_2,t_N)|... |\Omega_{t}(\eta_{N-1},t_1),...,\Omega_{t}(\eta_{N-1},t_N)]^T
                \end{align}
               Eq.~\ref{Eq.2.25} is the discretized form of Eqs.~\ref{Eq.2.15}-\ref{Eq.2.17}, using the Chebyshev pseudospectral method.\newline 
                Finally, we have
                         \begin{equation}\label{Eq.2.27}
                          G(V) = 0.
                          \end{equation}
                   The solution to the nonlinear equation system \ref{Eq.2.27} may be obtained by employing the Newton-Raphson technique 
                
                      \begin{eqnarray}
                        V_{n+1}=V_{n}-[J_{n}]^{-1}G(V_{n}),
                       \end{eqnarray}
                with the following Jacobian matrix
     \begin{align}
     \nonumber &  J =  
      \begin{bmatrix}
      \begin{vmatrix}  \frac{\partial G_{11}}{\partial V_{11}}\hspace{0.31cm}& \cdots \hspace{0.31cm}& \frac{\partial G_{11}}{\partial V_{1N}} \hspace{0.31cm}& \frac{\partial G_{11}}{\partial V_{21}} \hspace{0.31cm} & \cdots \hspace{0.31cm}& \frac{\partial G_{11}}{\partial V_{2N}}\hspace{0.31cm}& |\cdots|\hspace{0.31cm}& \frac{\partial G_{11}}{\partial V_{(N-1)1}}\hspace{0.31cm}& \cdots\hspace{0.31cm}& \frac{\partial G_{11}}{\partial V_{(N-1)N}} \\
     \ddots\hspace{0.31cm}& \ddots\hspace{0.31cm}& \ddots \hspace{0.31cm}& \ddots \hspace{0.31cm}&\ddots \hspace{0.31cm}& \ddots\hspace{0.31cm}&|\ddots|\hspace{0.31cm}&\ddots\hspace{0.31cm}&\ddots\hspace{0.31cm}&\ddots\\
             \frac{\partial G_{1N}}{\partial V_{11}}\hspace{0.31cm}& \cdots \hspace{0.31cm}& \frac{\partial G_{1N}}{\partial V_{1N}} \hspace{0.31cm}& \frac{\partial G_{1N}}{\partial V_{21}} \hspace{0.31cm}&\cdots \hspace{0.31cm}& \frac{\partial G_{1N}}{\partial V_{2N}}\hspace{0.31cm}&|\cdots|\hspace{0.31cm}&\frac{\partial G_{1N}}{\partial V_{(N-1)1}}\hspace{0.31cm}&...\hspace{0.31cm}&\frac{\partial G_{1N}}{\partial V_{(N-1)N}} 
      \end{vmatrix} \\
 \begin{vmatrix}   \frac{\partial G_{21}}{\partial V_{11}}\hspace{0.31cm}& \cdots \hspace{0.31cm}& \frac{\partial G_{21}}{\partial V_{1N}} \hspace{0.31cm}& \frac{\partial G_{21}}{\partial V_{21}} \hspace{0.31cm}&\cdots \hspace{0.31cm}& \frac{\partial G_{21}}{\partial V_{2N}}\hspace{0.31cm}&|\cdots|\hspace{0.31cm}&\frac{\partial G_{21}}{\partial V_{(N-1)1}}\hspace{0.31cm}&\cdots\hspace{0.31cm}&\frac{\partial G_{21}}{\partial V_{(N-1)N}} \\
\ddots\hspace{0.31cm}& \ddots\hspace{0.31cm}& \ddots \hspace{0.31cm}& \ddots \hspace{0.31cm}&\ddots \hspace{0.31cm}& \ddots\hspace{0.31cm}&|\ddots|\hspace{0.31cm}&\ddots\hspace{0.31cm}&\ddots\hspace{0.31cm}&\ddots\\
  \frac{\partial G_{2N}}{\partial V_{11}}\hspace{0.31cm}& \cdots \hspace{0.31cm}& \frac{\partial G_{2N}}{\partial V_{1N}} \hspace{0.31cm}& \frac{\partial G_{2N}}{\partial V_{21}} \hspace{0.31cm}&\cdots \hspace{0.31cm}& \frac{\partial G_{2N}}{\partial V_{2N}}\hspace{0.31cm}&|\cdots|\hspace{0.31cm}&\frac{\partial G_{2N}}{\partial V_{(N-1)1}}\hspace{0.31cm}&...\hspace{0.31cm}&\frac{\partial G_{2N}}{\partial V_{(N-1)N}} 
  \end{vmatrix} \\ 
 \begin{vmatrix} \hspace{0.5cm} \vdots& \hspace{0.95cm}\vdots & \hspace{0.95cm} \vdots &\hspace{0.95cm} \vdots &\hspace{0.95cm} \vdots &\hspace{0.95cm} \vdots&\hspace{0.95cm} \vdots&\hspace{0.95cm} \vdots&\hspace{0.95cm} \vdots&\hspace{0.95cm} \vdots \hspace{0.5cm} \end{vmatrix} \\ 
 \begin{vmatrix}  \frac{\partial G_{(N-1)1}}{\partial V_{11}}& \cdots & \frac{\partial G_{(N-1)1}}{\partial V_{1N}} & \frac{\partial G_{(N-1)1}}{\partial V_{21}} &\cdots & \frac{\partial G_{(N-1)1}}{\partial V_{2N}}&|\cdots|&\frac{\partial G_{(N-1)1}}{\partial V_{(N-1)1}}&\cdots&\frac{\partial G_{(N-1)1}}{\partial V_{(N-1)N}}  \\ 
\ddots& \ddots& \ddots & \ddots &\ddots & \ddots&|\ddots|&\ddots&\ddots&\ddots\\
  \frac{\partial G_{(N-1)N}}{\partial V_{11}}& \cdots & \frac{\partial G_{(N-1)N}}{\partial V_{1N}} & \frac{\partial G_{(N-1)N}}{\partial V_{21}} &\cdots & \frac{\partial G_{(N-1)N}}{\partial V_{2N}}&|\cdots|&\frac{\partial G_{(N-1)N}}{\partial V_{(N-1)1}}&...&\frac{\partial G_{(N-1)N}}{\partial V_{(N-1)N}} \end{vmatrix} \\ 
      \end{bmatrix}
    \end{align} 
 where the order of the Jacobian matrix is $(N-1)N \times (N-1)N$.

 \section{Stability analysis of the gBF equation} \label{sec3}
 In order to conduct a study of the method's stability for the gBF equation, we make use of the relation between the collocation points and the Gauss-type integration formulas. These quadrature formulae have the characteristics that let us make the transition from summing to integrating.
The gBF equation is given by
\begin{align} \label{Eq.3.1}
U_t(\eta,t) + A(U^{\delta + 1}(\eta,t))_{\eta} - B U_{\eta \eta}(\eta,t) - C F(U(\eta,t)) = 0.
\end{align}
Multiplying Eq.~\ref{Eq.3.1} by $U(\eta_\alpha,t)$ with the weight function $w(\eta_\alpha) = \frac{1}{\sqrt(1-\eta_{\alpha}^2)}$ and summing, we have 
\begin{align}\label{Eq.3.2}
\sum_{\alpha=0}^{N} U(\eta_\alpha, t) w(\eta_\alpha) (U_t(\eta_\alpha,t) + A(U^{\delta + 1}(\eta_\alpha,t))_\eta - B U_{\eta \eta}(\eta_\alpha,t)- C F(U(\eta_\alpha,t))) = 0,
\end{align}
where $A$, $B$, and $C$ are constants. Eq.~\ref{Eq.3.2} can be written as
\begin{align}\label{Eq.3.3}
\nonumber \frac{1}{2}\frac{d}{d t}\sum_{\alpha=0}^{N} U^2(\eta_\alpha, t) w(\eta_\alpha) &= -A \sum_{\alpha=0}^{N} U(\eta_\alpha, t) (U^{\delta + 1}(\eta_\alpha,t))_\eta w(\eta_\alpha)\\
& + B \sum_{\alpha=0}^{N} U(\eta_\alpha, t) U_{\eta \eta}(\eta_\alpha,t) w(\eta_\alpha)+ C \sum_{\alpha=0}^{N} U(\eta_\alpha, t) F(U(\eta_\alpha,t)) w(\eta_\alpha).
\end{align}
The first two terms of the right hand side in Eq. \ref{Eq.3.3} can be shown easily non-positive using integration by parts as follows:\newline 
Let $w U = z$ $\Rightarrow$ $U_{\eta} = \frac{(wU)_{\eta}}{w} - \frac{w_{\eta}U}{w} = \frac{z_{\eta}}{w} + z(\frac{1}{w})_{\eta}$, which leads to the following result  
\begin{align}\label{Eq.3.7}
\nonumber \sum_{\alpha=0}^{N}   U(\eta_\alpha, t) U_{\eta \eta}(\eta_\alpha,t) w(\eta_\alpha) & = \int_{-1}^{1}  w(\eta) U(\eta, t) U_{\eta \eta}(\eta,t) d\eta\\
& \nonumber = -\int_{-1}^{1} z_{\eta} (\frac{z_{\eta}}{w} + z(\frac{1}{w})_{\eta}) d \eta \\
& = -\int_{-1}^{1} \frac{z_{\eta}^2}{w} d \eta + \frac{1}{2} \int_{-1}^{1} z^2 (\frac{1}{w})_{\eta \eta} d \eta ~\leq 0.
\end{align}
Hence the Eq.~\ref{Eq.3.3} becomes
\begin{align}\label{Eq.3.8}
 \frac{1}{2}\frac{d}{d t}\sum_{\alpha=0}^{N} U^2(\eta_\alpha, t) w(\eta_\alpha)  \leq C \sum_{\alpha=0}^{N} U(\eta_\alpha, t) F(U(\eta_\alpha,t)) w(\eta_\alpha).
\end{align}
Since, we have $F(U(\eta,t)) = U(\eta,t)(1-U^{\delta}(\eta,t))$ and hence the above Eq.~\ref{Eq.3.8} transformed into the following equation 
\begin{align*}
\frac{1}{2}\frac{d}{d t}\sum_{\alpha=0}^{N} U^2(\eta_\alpha, t) w(\eta_\alpha) \leq C \sum_{\alpha=0}^{N} U^2(\eta_\alpha, t) w(\eta_\alpha) -C \sum_{\alpha=0}^{N} U^{\delta + 2}(\eta_\alpha,t) w(\eta_\alpha),
\end{align*}
\begin{align*}
\Rightarrow \frac{\frac{d}{d t}\sum_{\alpha=0}^{N} U^2(\eta_\alpha, t) w(\eta_\alpha)}{\sum_{\alpha=0}^{N} U^2(\eta_\alpha, t) w(\eta_\alpha)} \leq 2C  - 2C \frac{\sum_{\alpha=0}^{N} U^{\delta + 2}(\eta_\alpha,t) w(\eta_\alpha)}{\sum_{\alpha=0}^{N} U^2(\eta_\alpha, t) w(\eta_\alpha)} \leq 2C,
\end{align*}
\begin{align*}
\Rightarrow log(\sum_{\alpha=0}^{N} U^2(\eta_\alpha, t) w(\eta_\alpha))  \leq 2C t + log(k).
\end{align*}
Further simplifying, we get 
\begin{align*}
\sum_{\alpha=0}^{N} U^2(\eta_\alpha, t) w(\eta_\alpha) \leq k \exp(2C t).
\end{align*}
Finally using the initial condition, we obtain
\begin{align*}
\sum_{\alpha=0}^{N} U^2(\eta_\alpha, t) w(\eta_\alpha) \leq \exp(2C t) \sum_{\alpha=0}^{N} U^2(\eta_\alpha, t_0) w(\eta_\alpha),
\end{align*}
hence conclude the result.

     \section{Numerical findings and discussion} \label{sec4}
    In order to illustrate the performance, we consider Eq.~\ref{Eq.1.4} with the finite spatial and time intervals with the following initial condition
    \begin{equation}
    U(\tilde{\eta},0) = \left(\frac{1}{2} + \frac{1}{2}  \tanh(U_1 \tilde{\eta})\right)^{\frac{1}{\delta}}, ~~~~0 \leq \tilde{\eta} \leq 1,
    \end{equation}    
    and boundary conditions
    \begin{align} 
   &U(0,\tilde{t})=\left(\frac{1}{2} + \frac{1}{2}  \tanh(U_1  U_2  \tilde{t})\right)^{\frac{1}{\delta}},~~~~ \tilde{t} \geq 0, \\
   & U(1,\tilde{t}) = \left(\frac{1}{2} + \frac{1}{2}  \tanh(U_1  (1 - U_2  \tilde{t}))\right)^{\frac{1}{\delta}},~~~~ \tilde{t} \geq 0.
    \end{align}
    The exact solution~\cite{wang1990exact} of Eq.~\ref{Eq.1.4} is
    \begin{align}
    U(\tilde{\eta},\tilde{t}) = \left(\frac{1}{2} + \frac{1}{2}  \tanh\left(U_1 (\tilde{\eta} - U_2  \tilde{t} )\right) \right)^{\frac{1}{\delta}} ,~~~~0 \leq \tilde{\eta} \leq 1,~~~~ \tilde{t} \geq 0,
    \end{align} where $U_1$ and $U_2$ are given by
    $$U_1= \frac{-\sigma_1  \delta}  {2 (\delta + 1)},~~~~U_2= \left(\frac{\sigma_1 }{\delta + 1} + \sigma_2   \frac{\delta + 1}{\sigma_1 }\right).$$
The accuracy and efficiency of the technique have been verified for the distinct values of constants $\sigma_1$ and $\sigma_2$ with different sets of time and various values of $\delta$, all of which yielded excellent and encouraging numerical results.
To express the error, $L_\infty$-norm has been calculated at different grid points which can be defined as 
     $$L_{\infty} = \text{max}\{|U_{i,j}^{exa} - U_{i,j}^{num}|: i=1,2,...,N-1,~j=1,2,...,N \},$$
     where $U^{exa}$ and $U^{num}$ denotes the exact and approximated solutions respectively.
     \subsection{Numerical examples for the gBF equation}	
\begin{figure}[htbp] 
        \centering
       \includegraphics[width=6.4cm, height=6.5cm]{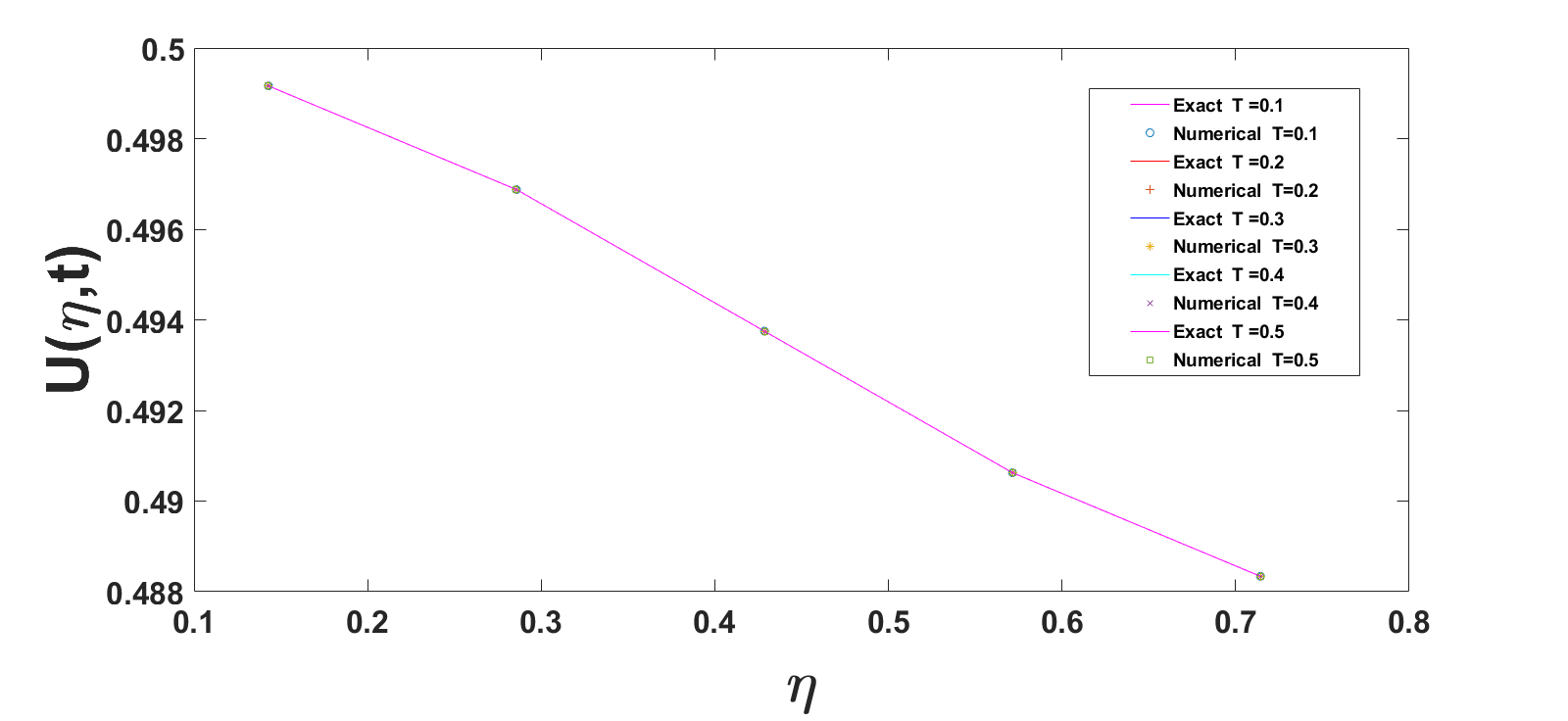} 
      \includegraphics[width=6.4cm, height=6.5cm]{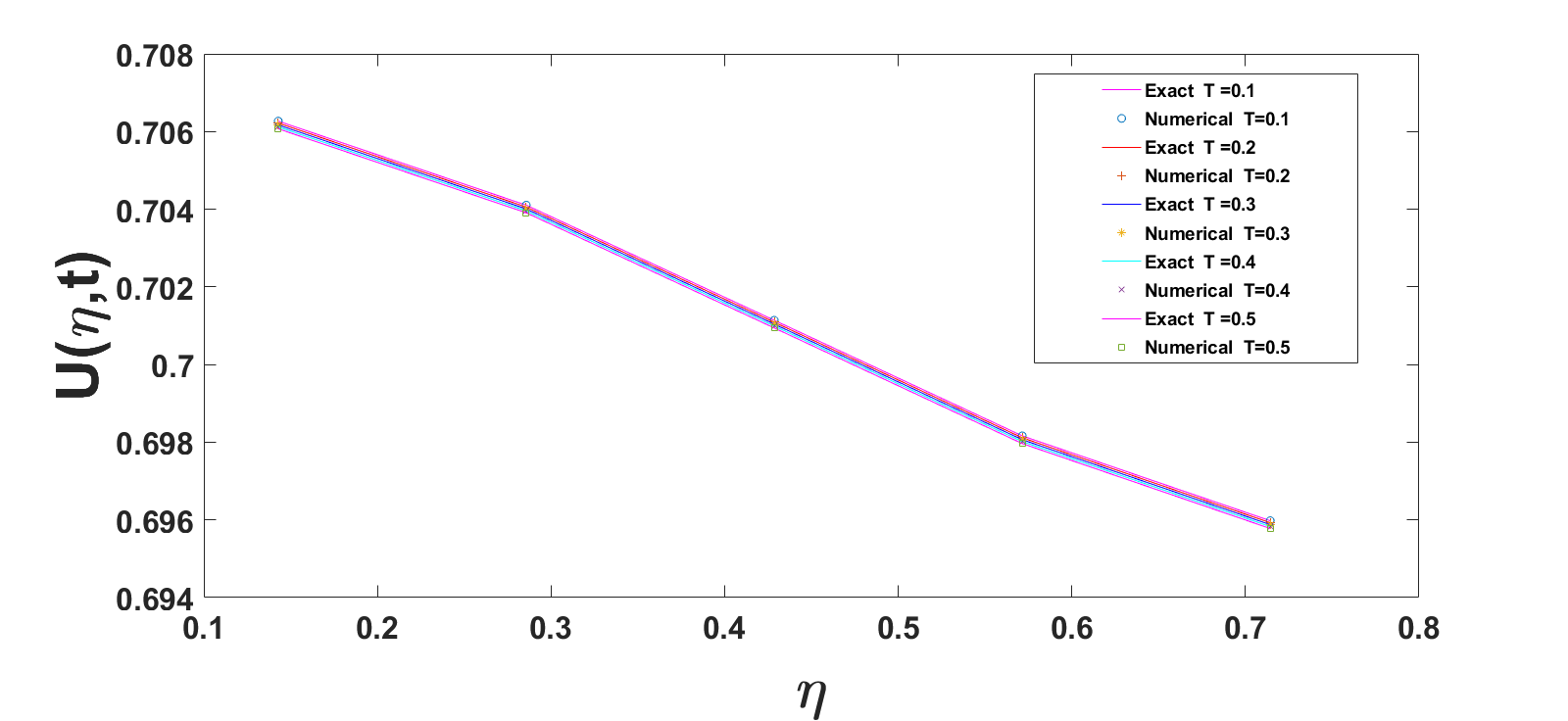}\\
      \includegraphics[width=6.4cm, height=6.5cm]{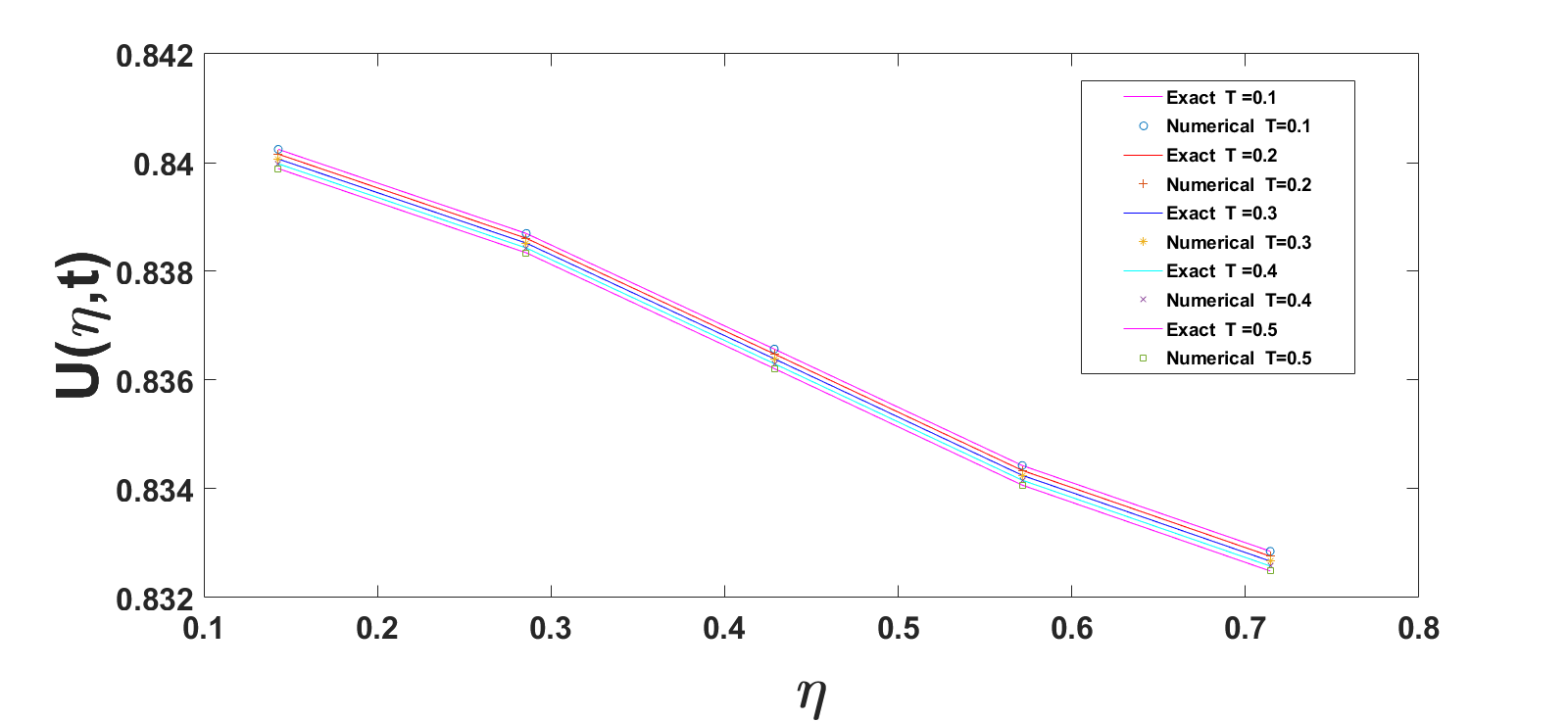}
       \includegraphics[width=6.4cm, height=6.5cm]{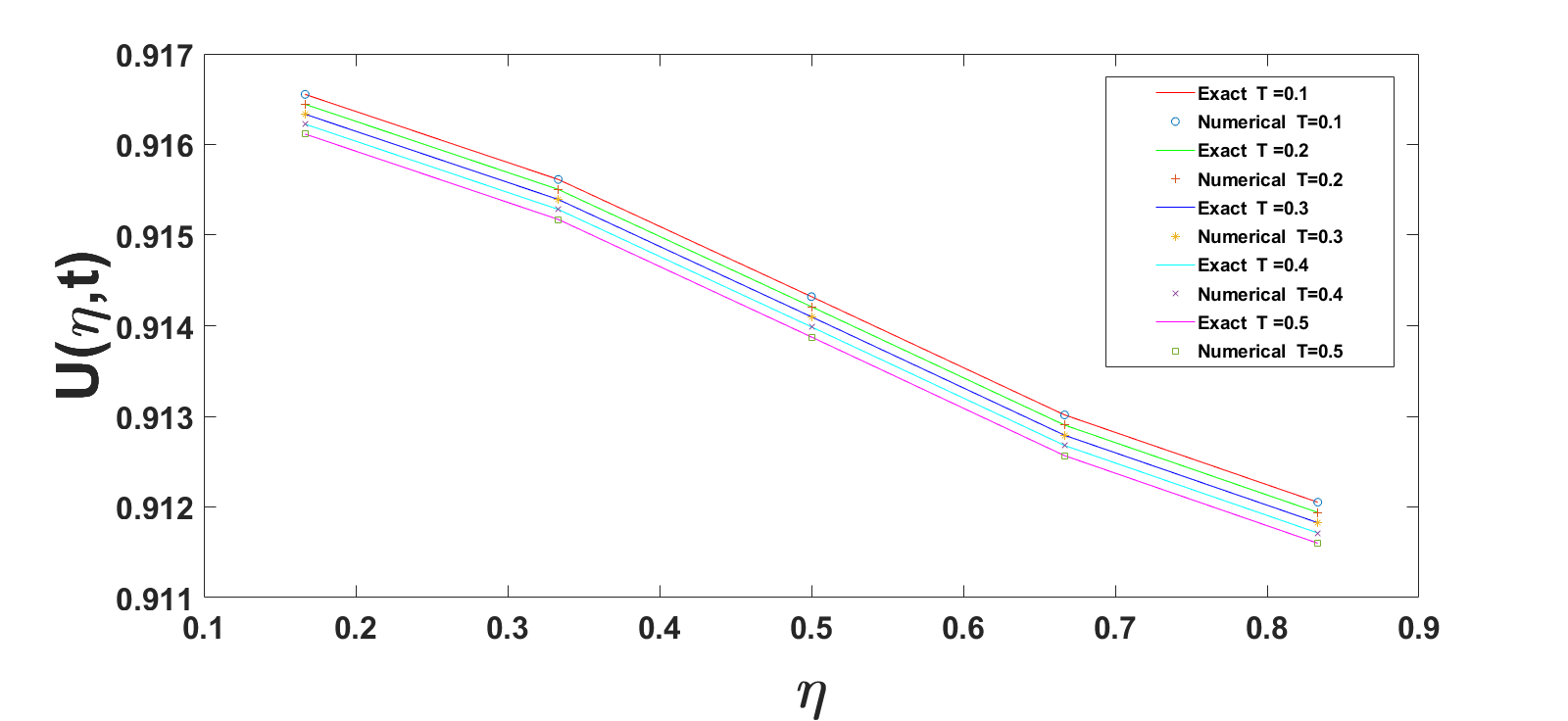}
       \caption{Numerical  vs Exact solution of $ U(\eta,t) $ for Example~\ref{eg1} when  different time $ T=$ 0.1, 0.2, 0.3, 0.4 and 0.5   for $\delta =$ 1, 2, 4, and 8, respectively.}
        \label{fig1_eg1}
       \end{figure}

\begin{figure}[htbp]
    \centering
     \includegraphics[width=14cm, height=5.3cm]{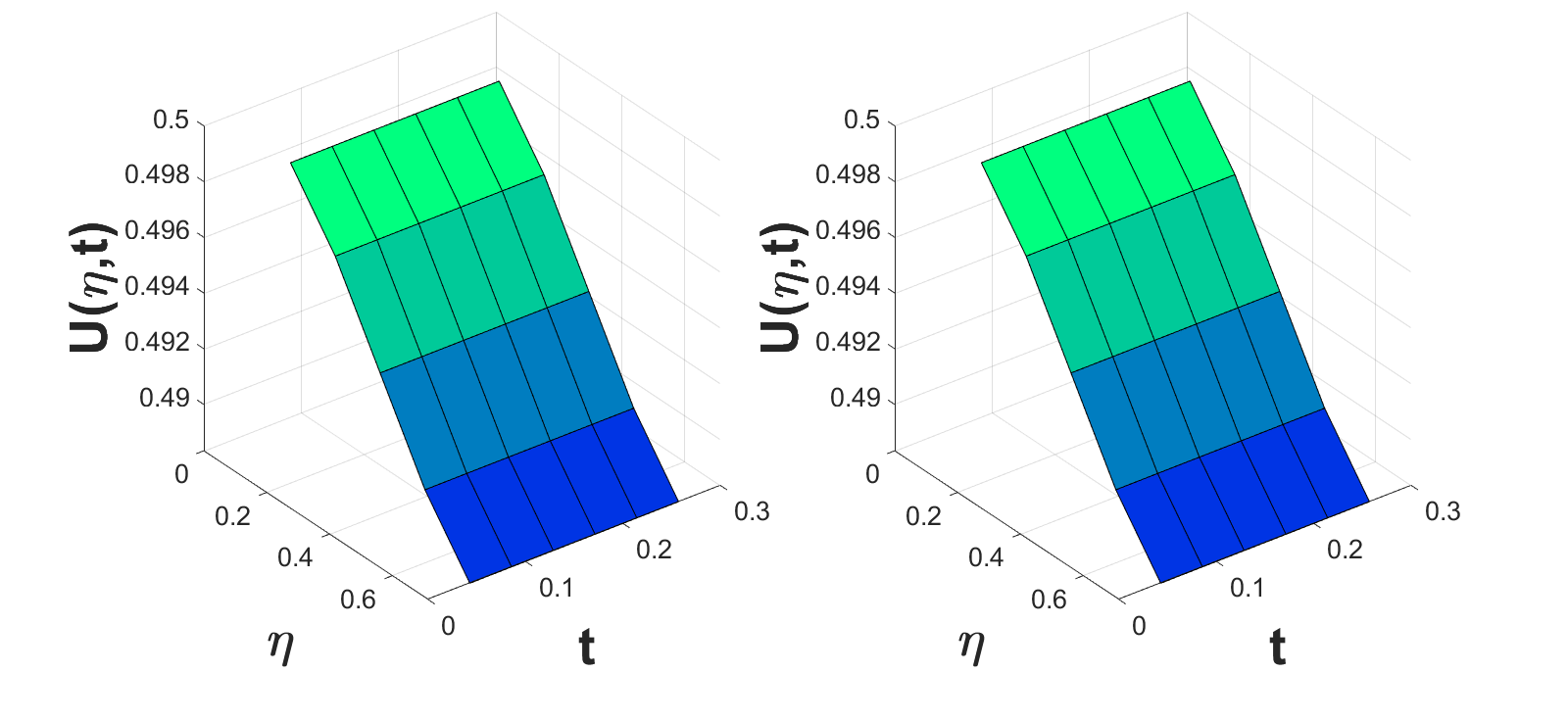} \\
     \includegraphics[width=14cm, height=5.3cm]{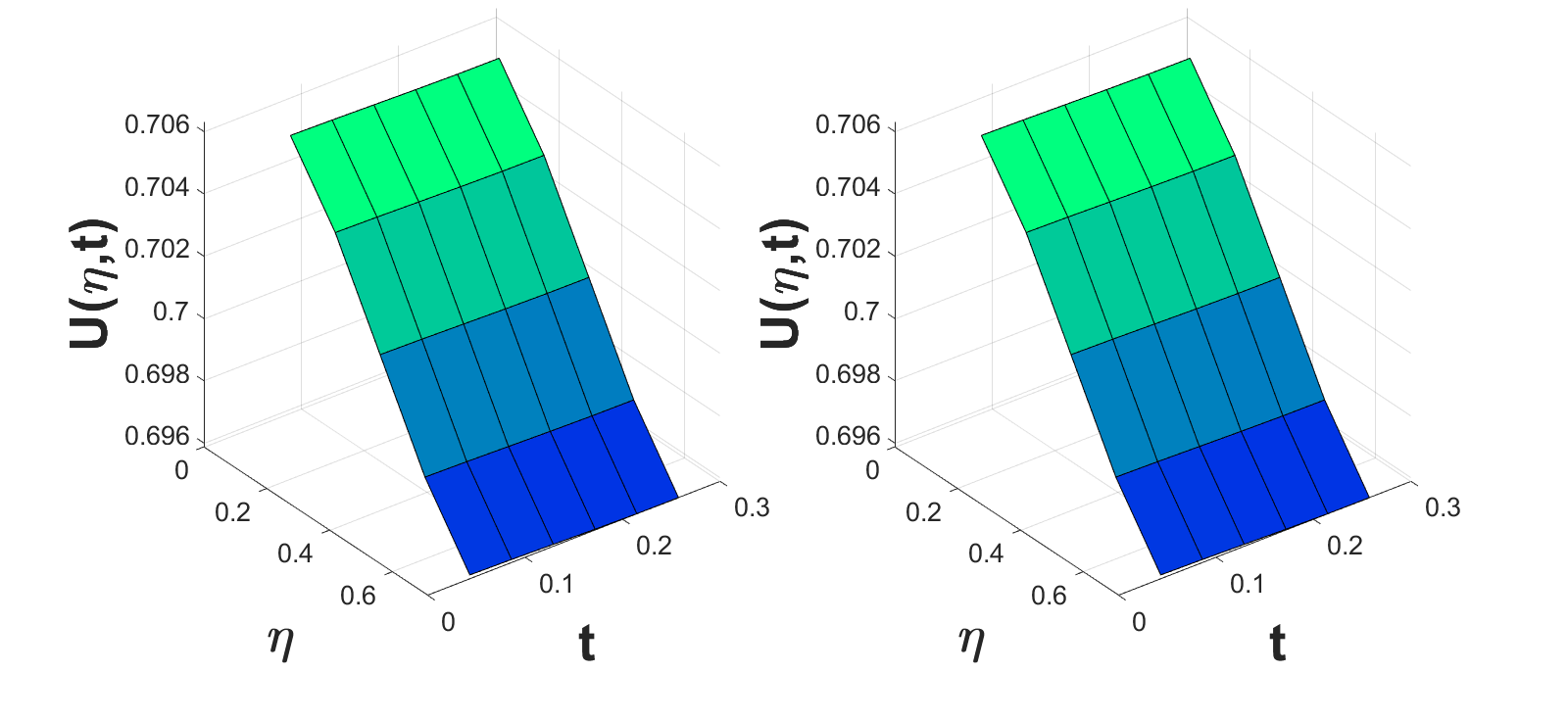}\\ 
      \includegraphics[width=14cm, height=5.3cm]{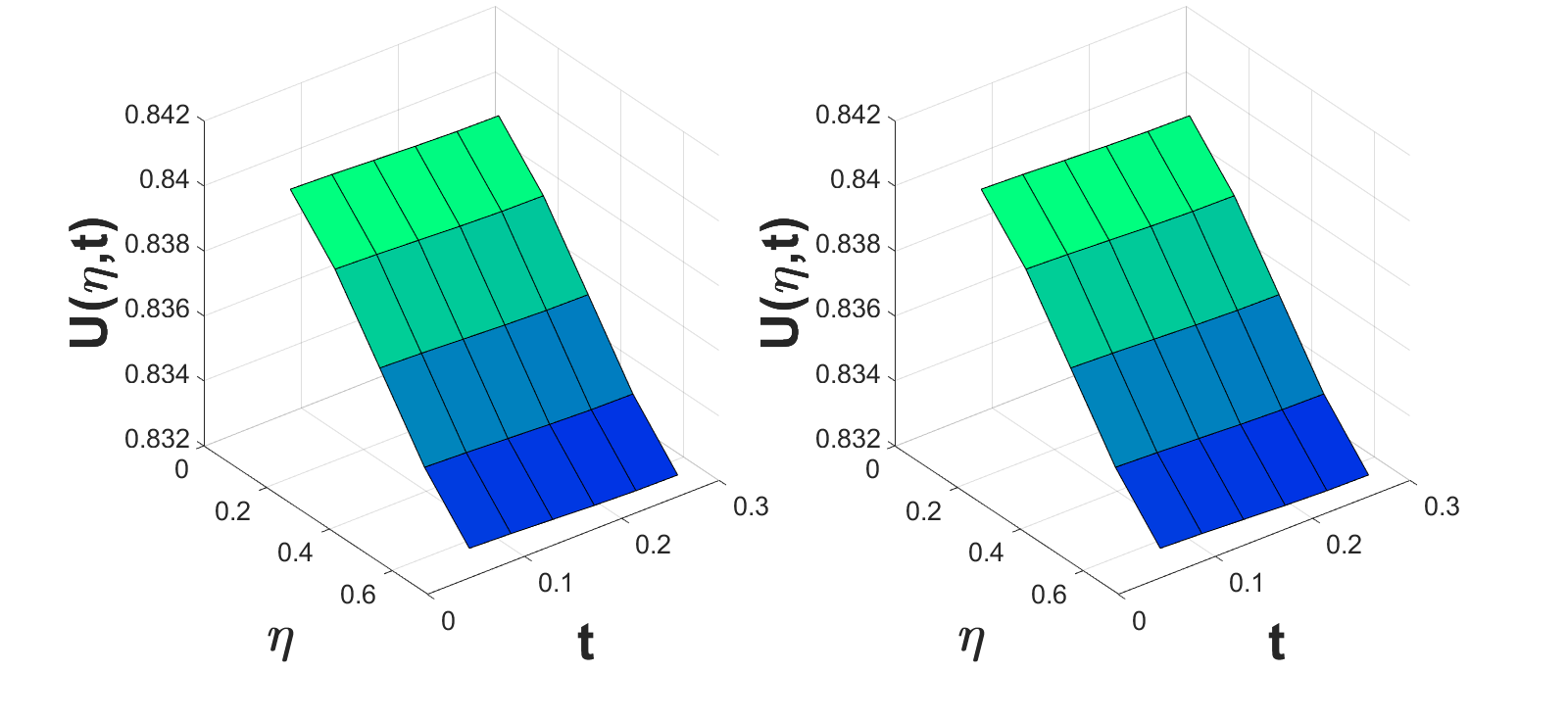}\\
      \includegraphics[width=14cm, height=5.3cm]{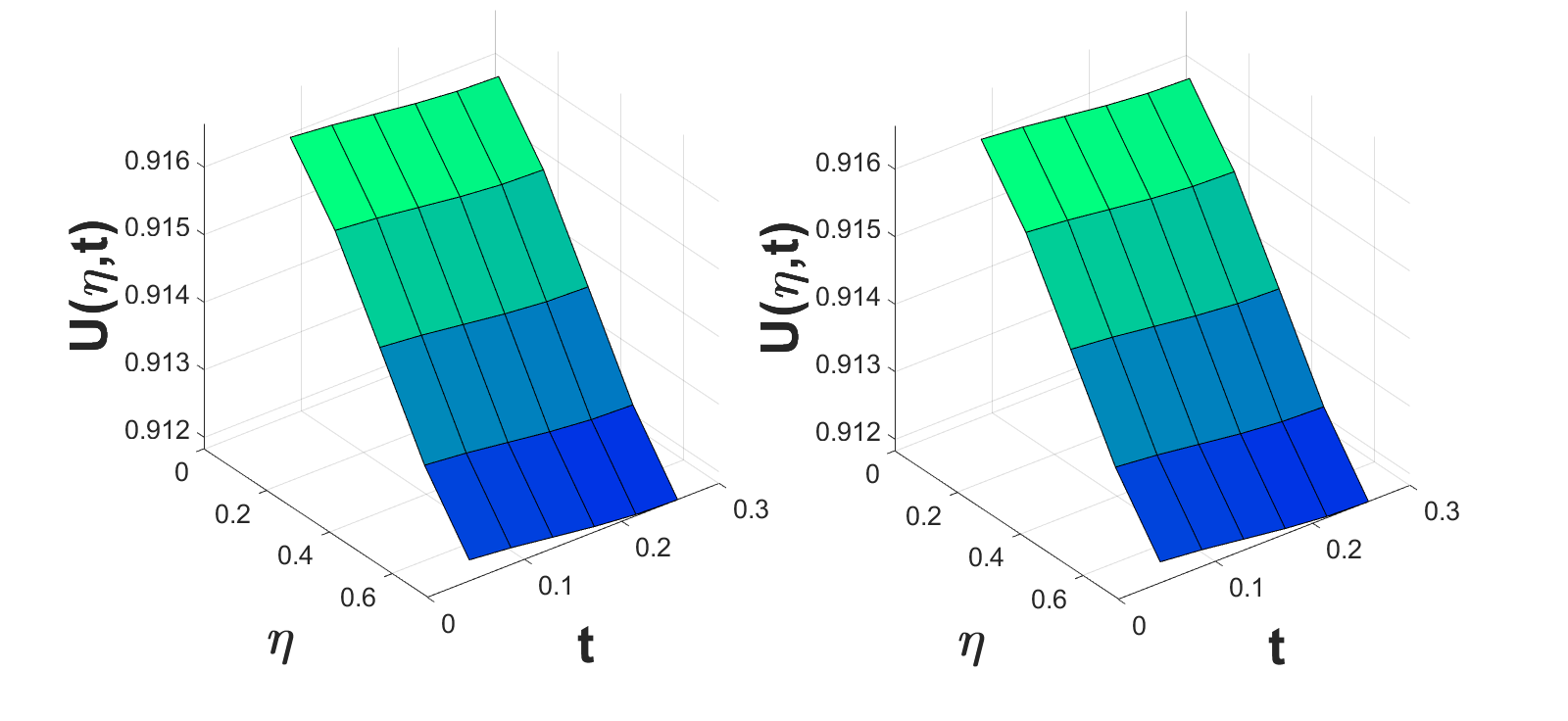}
    \caption{3D graphs for exact(left) vs approximate(right) solution for Example~\ref{eg1} at time $ T=$ 0.3 for $\delta =$ 1, 2, 4, and 8, respectively.}
    \label{fig2_eg1}
  \end{figure}

\begin{figure}[h!] 
           \centering
          \includegraphics[width=6.4cm, height=6.5cm]{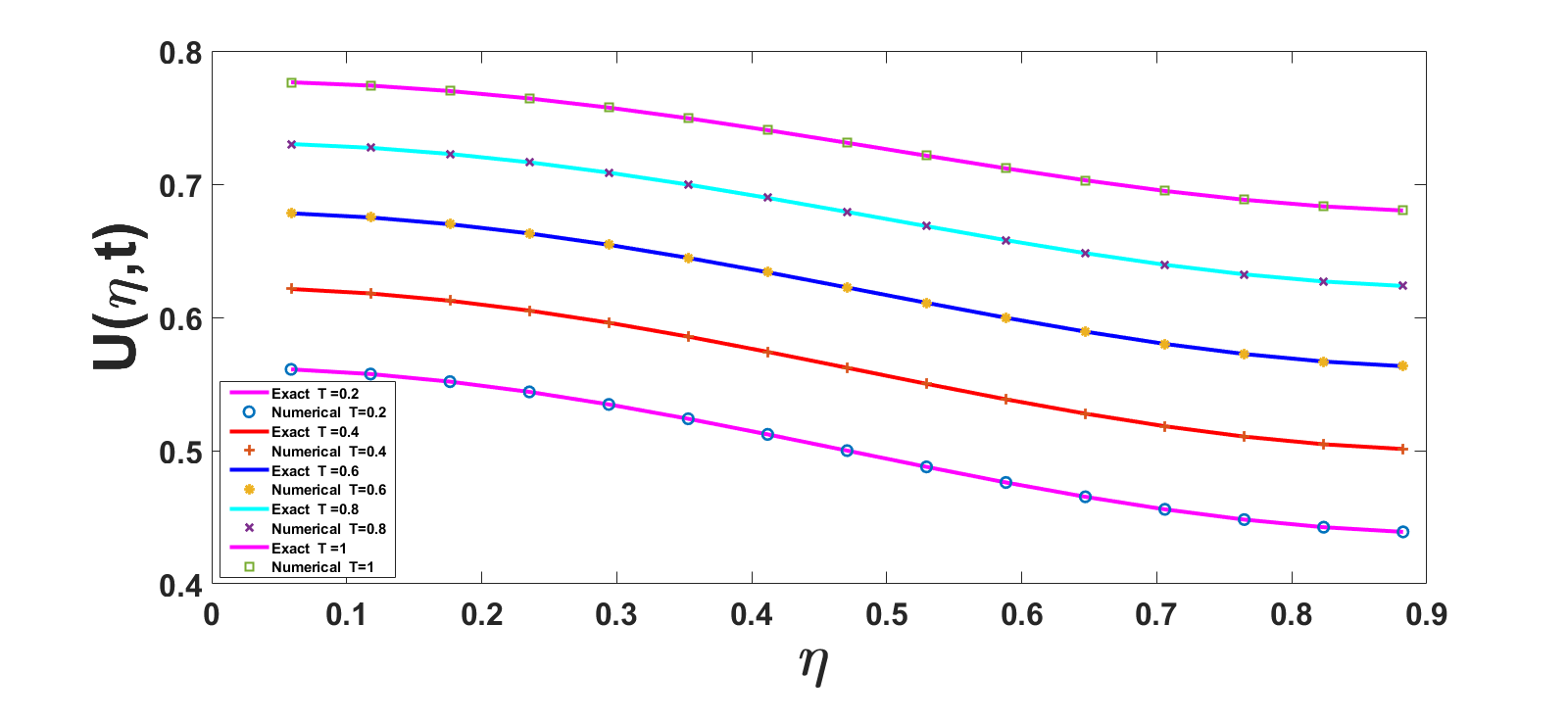} 
         \includegraphics[width=6.4cm, height=6.5cm]{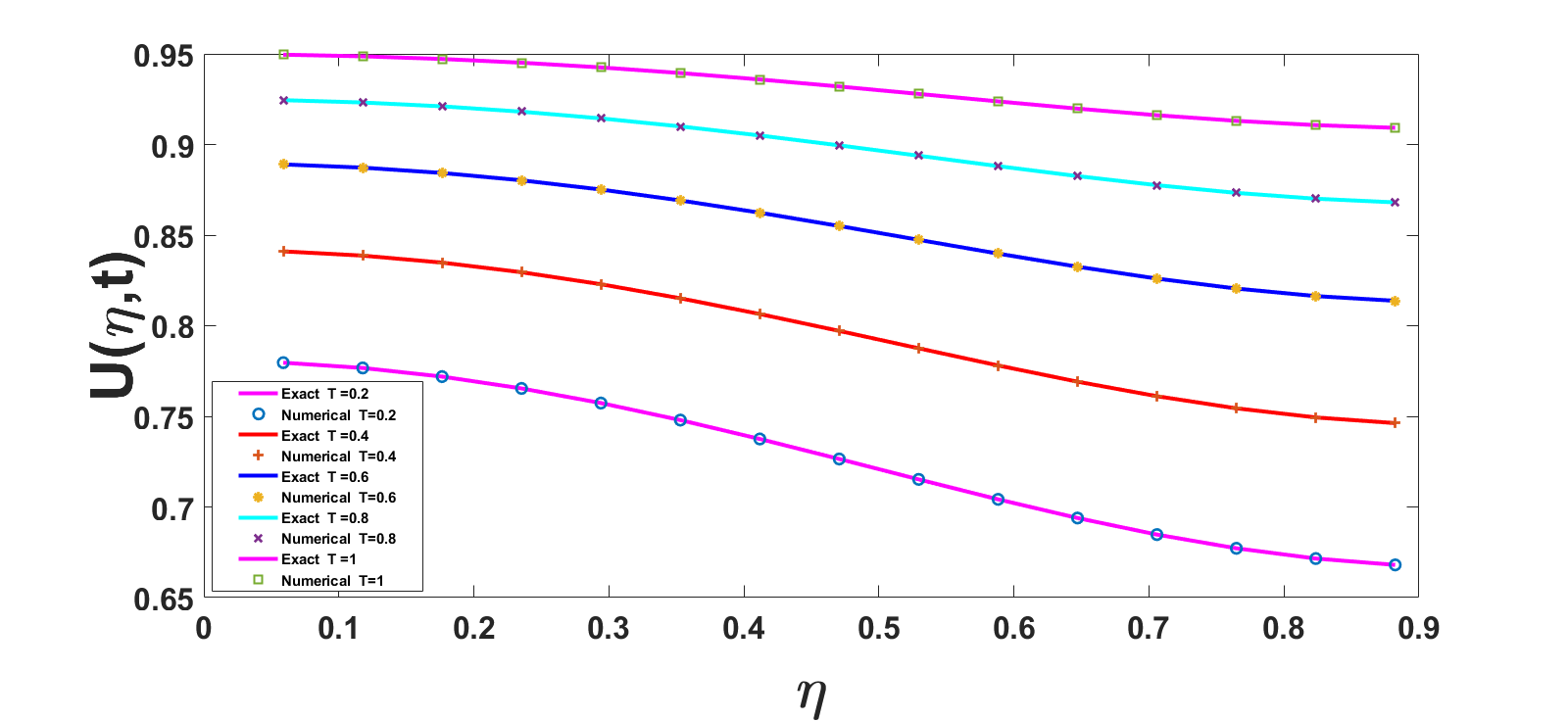}\\
         \includegraphics[width=6.4cm, height=6.5cm]{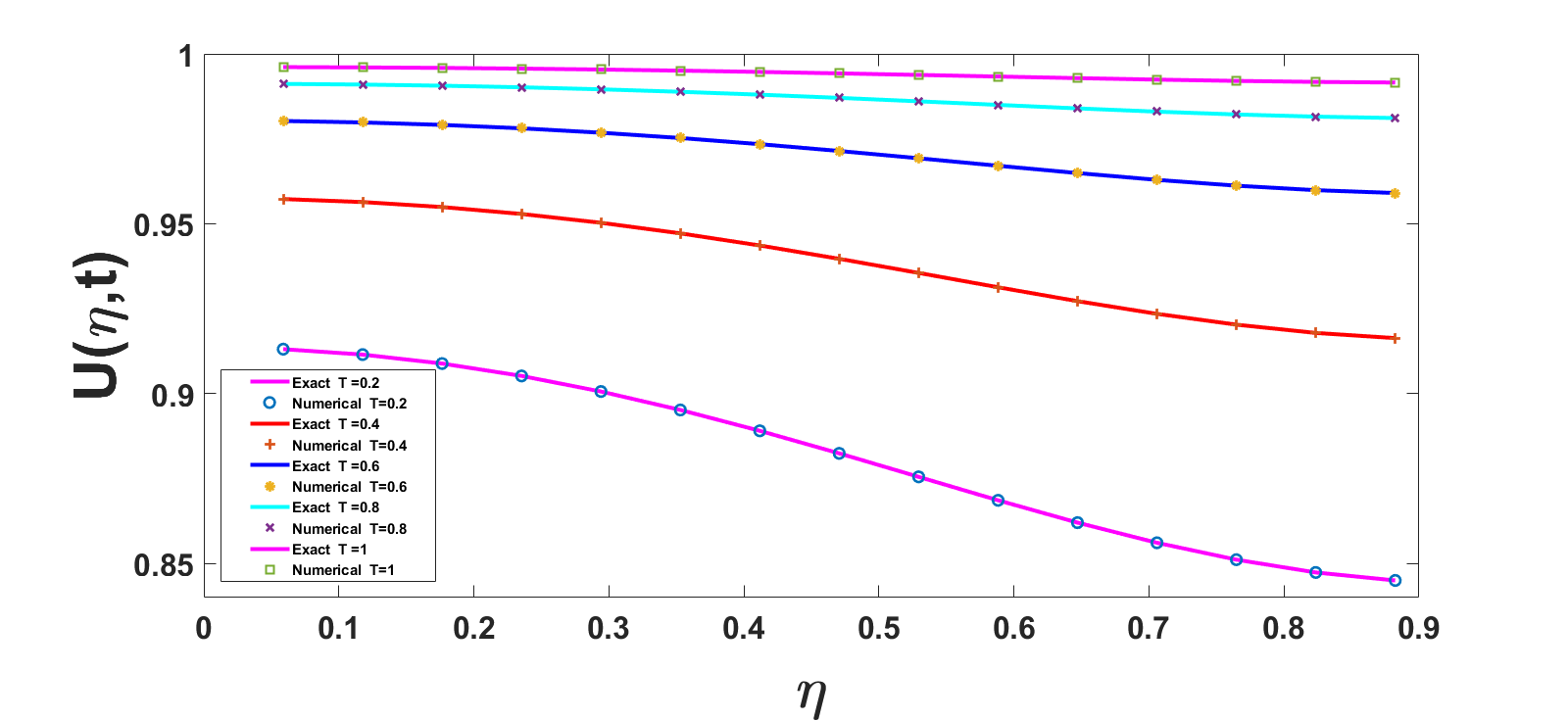}
          \includegraphics[width=6.4cm, height=6.5cm]{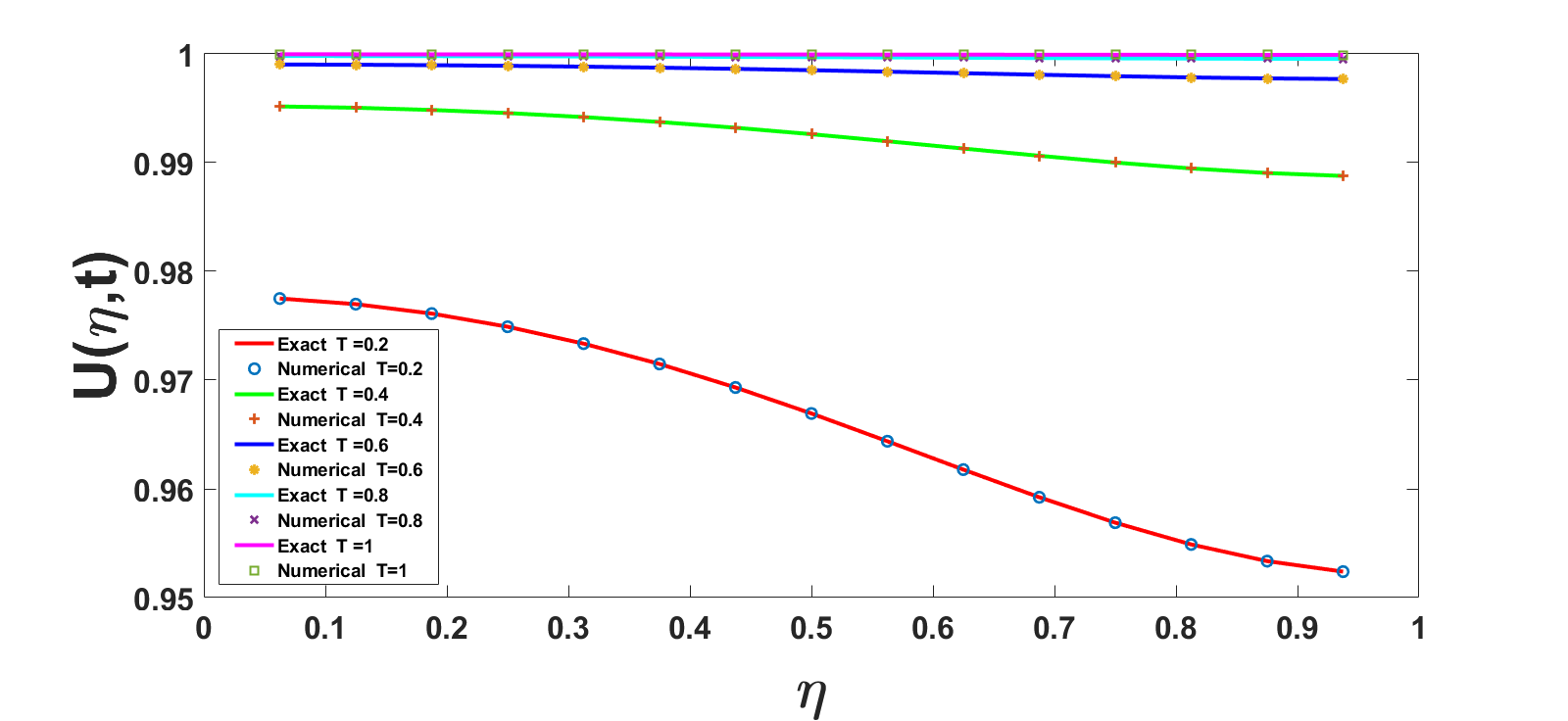}
          \caption{Numerical  vs Exact solution of $ U(\eta,t) $ for Example~\ref{eg2} when  different time $ T=$ 0.2, 0.4, 0.6, 0.8 and 1   for $\delta =$ 1, 2, 4, and 8, respectively.}
           \label{fig3a_eg2}
          \end{figure}

\begin{figure}[htbp]
    \centering
     \includegraphics[width=14cm, height=5.3cm]{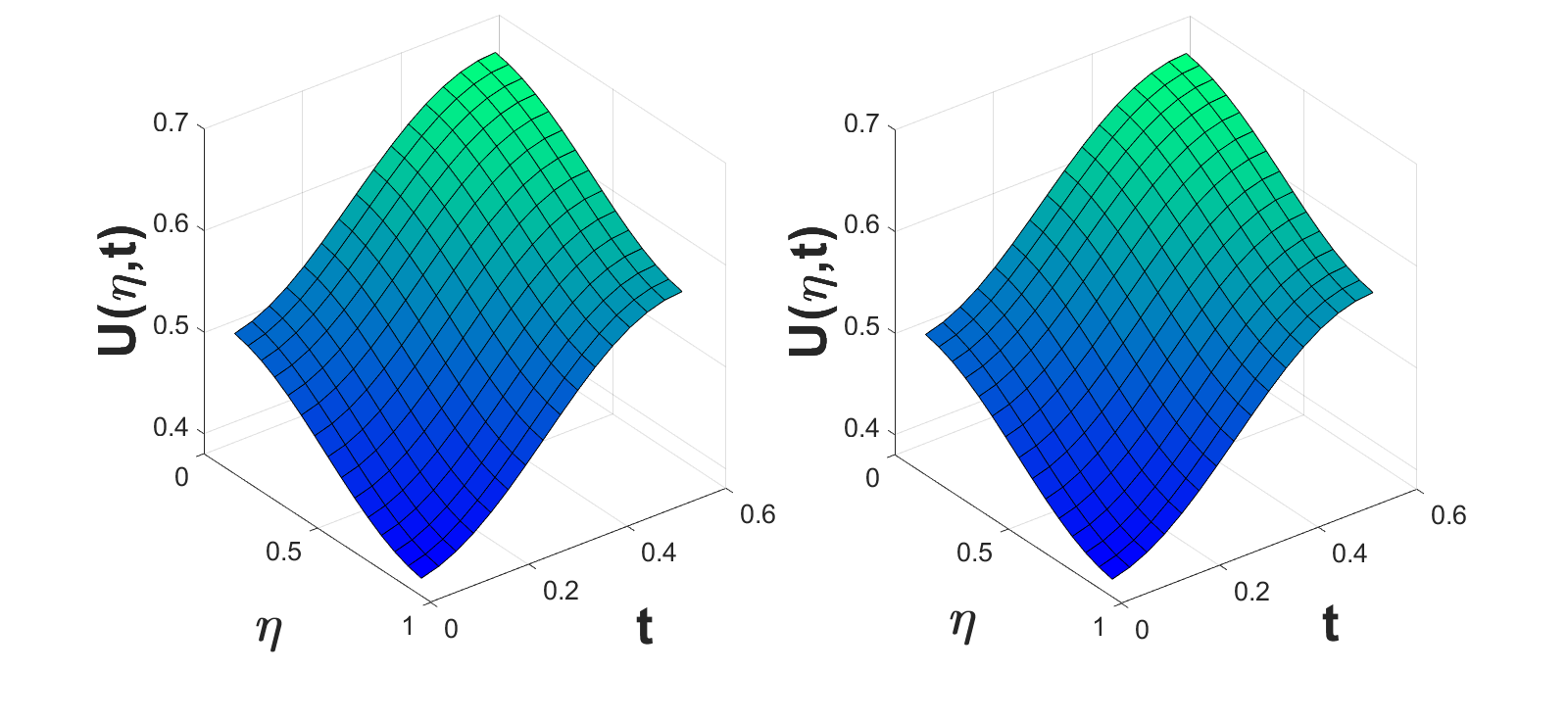} \\
     \includegraphics[width=14cm, height=5.3cm]{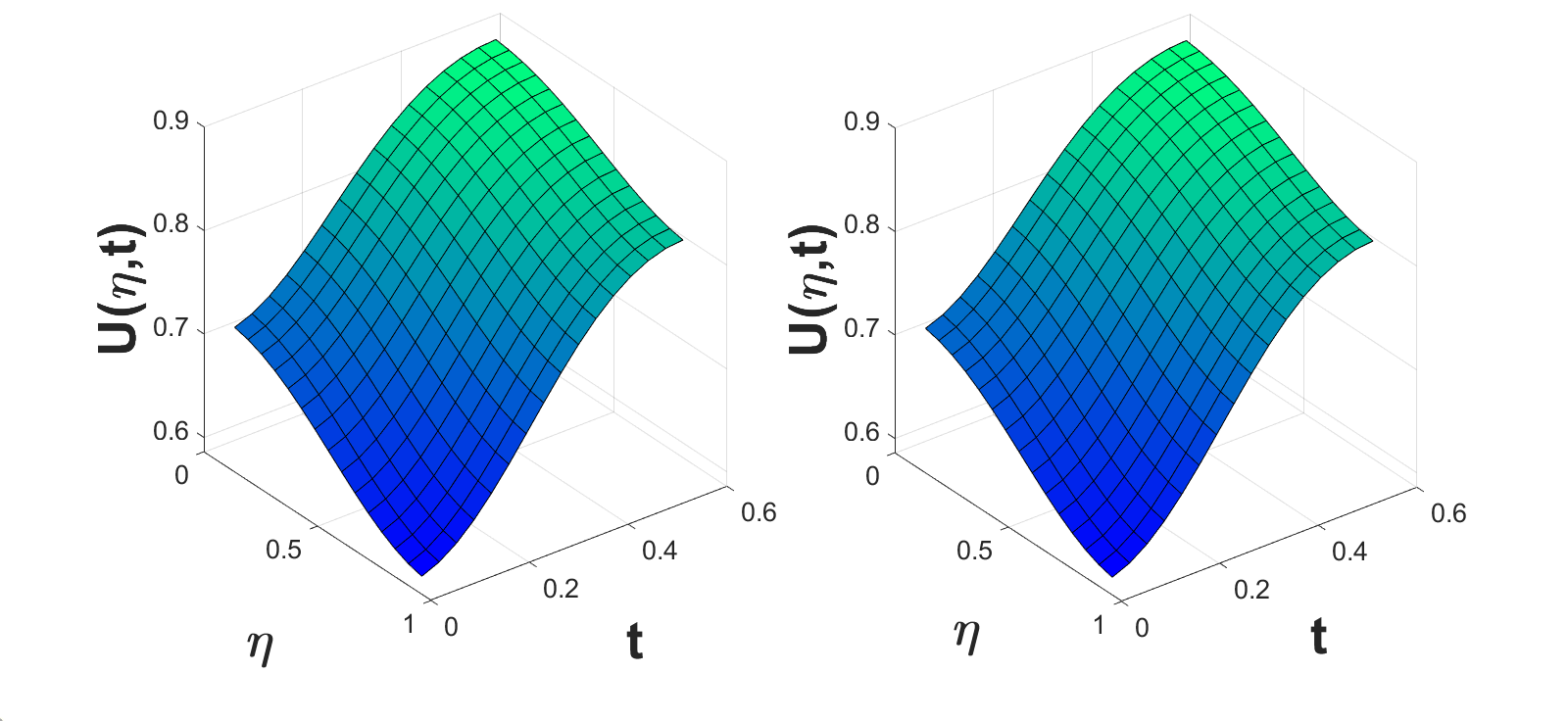}\\ 
      \includegraphics[width=14cm, height=5.3cm]{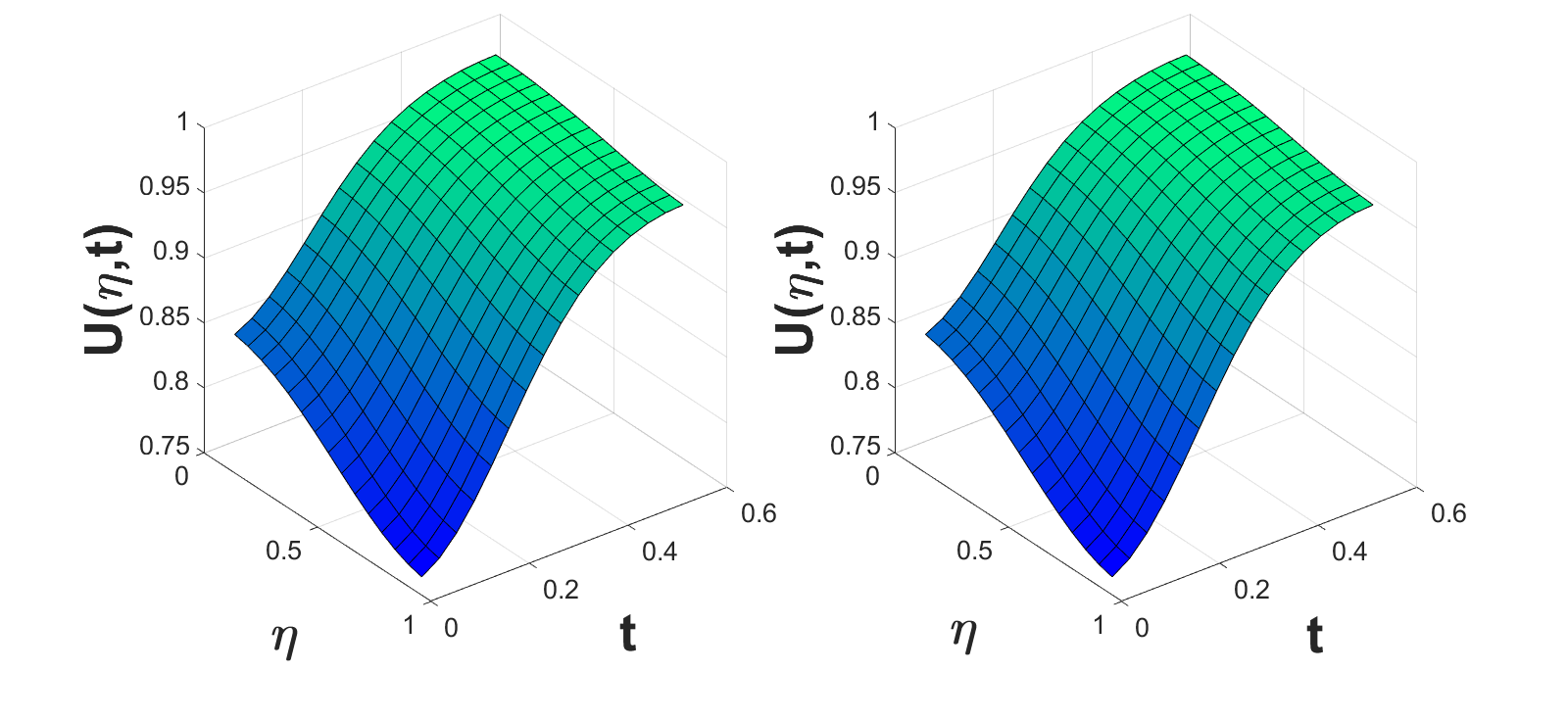}\\
      \includegraphics[width=14cm, height=5.3cm]{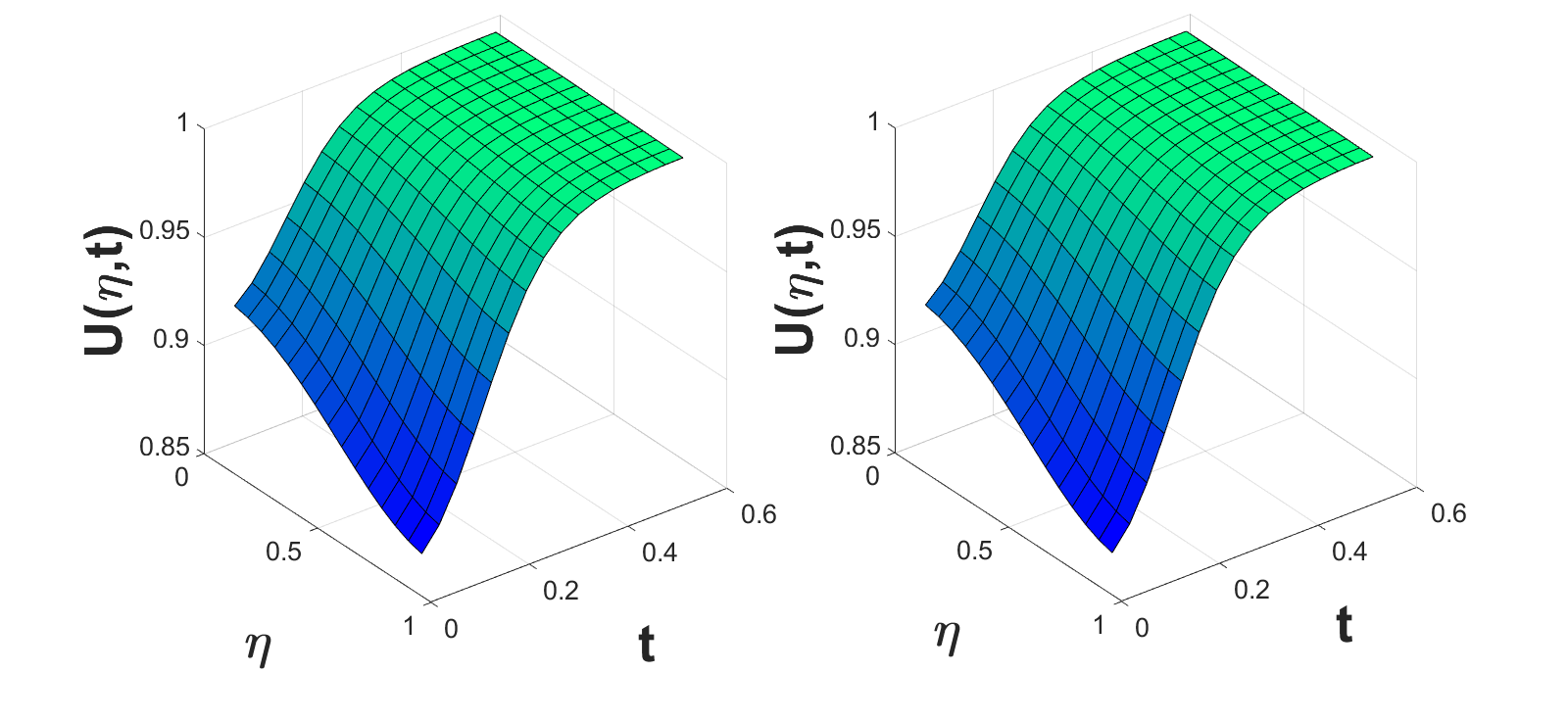}
    \caption{3D graphs for exact(left) vs approximate(right) solution for Example~\ref{eg2} at time $ T=$ 0.6 for $\delta =$ 1, 2, 4, and 8, respectively.}
    \label{fig3b_eg2}
  \end{figure}

\begin{figure}[htbp]
         \centering
        \includegraphics[width=0.93\linewidth]{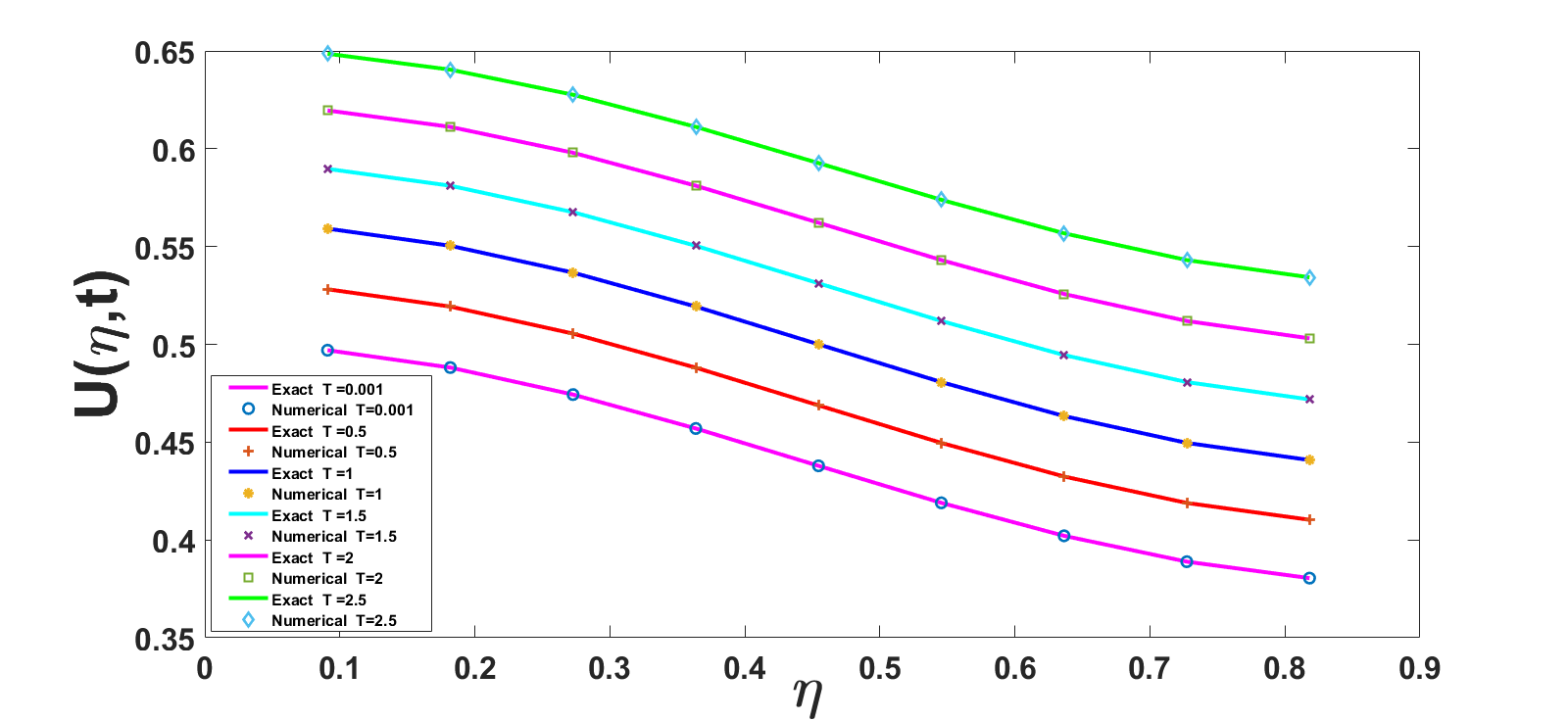} 
       \includegraphics[width=0.93\linewidth]{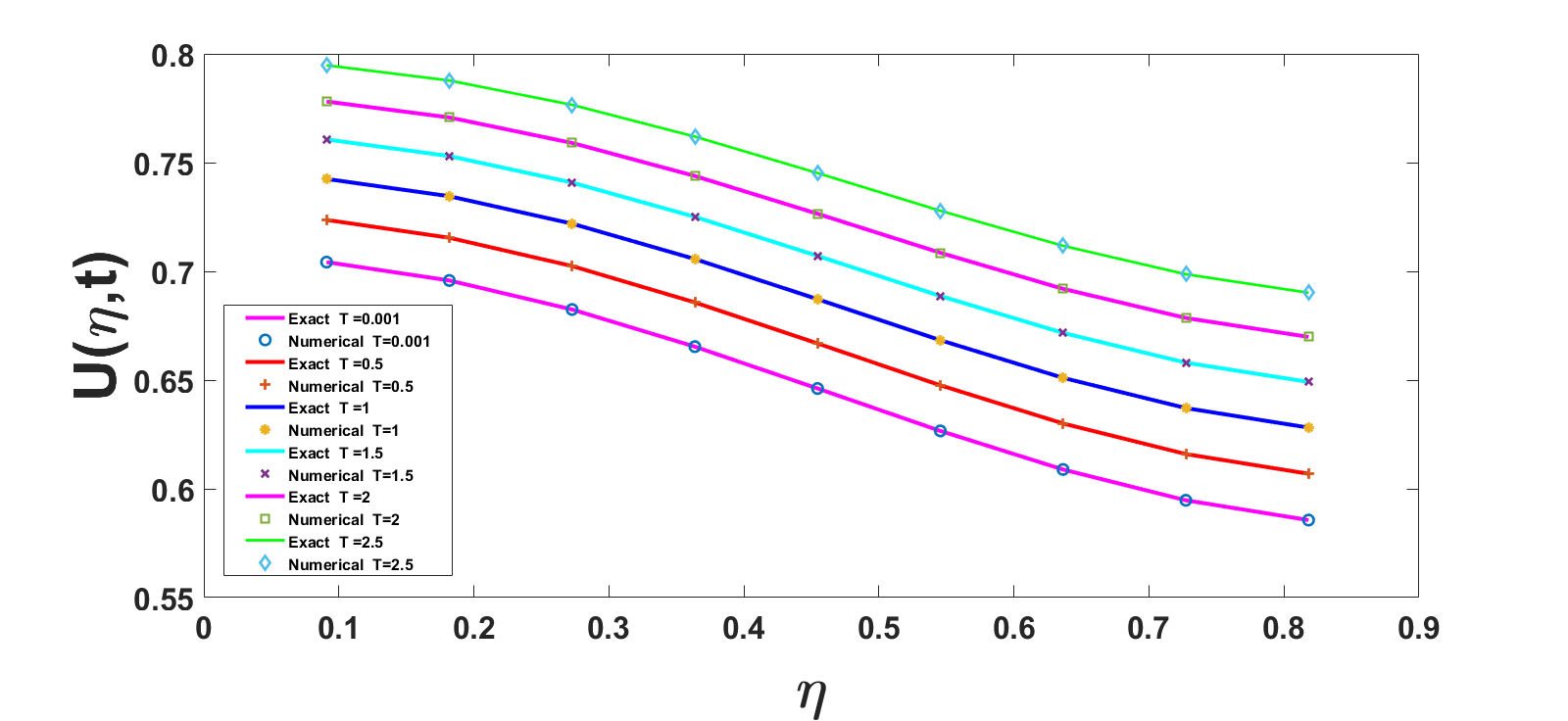}
       \includegraphics[width=0.93\linewidth]{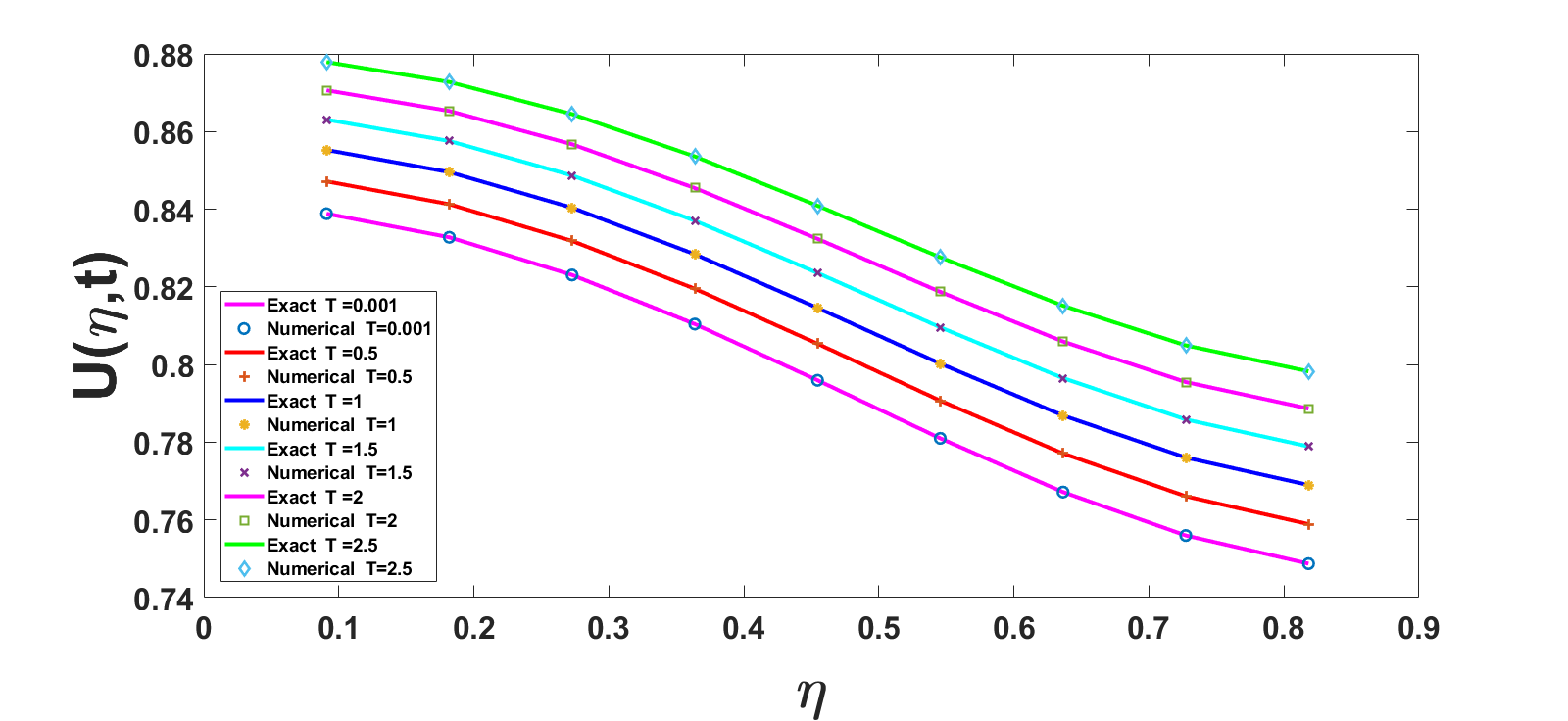}
        \caption{Numerical  vs Exact solution for Example~\ref{eg3} when  different time $ T=$ 0.001, 0.5, 1, 1.5, 2 and 2.5   for $\delta =$ 1, 2, and 4, respectively.}
         \label{fig4_eg3}
        \end{figure}

\begin{figure}[htbp]
    \centering
     \includegraphics[width=14cm, height=6.3cm]{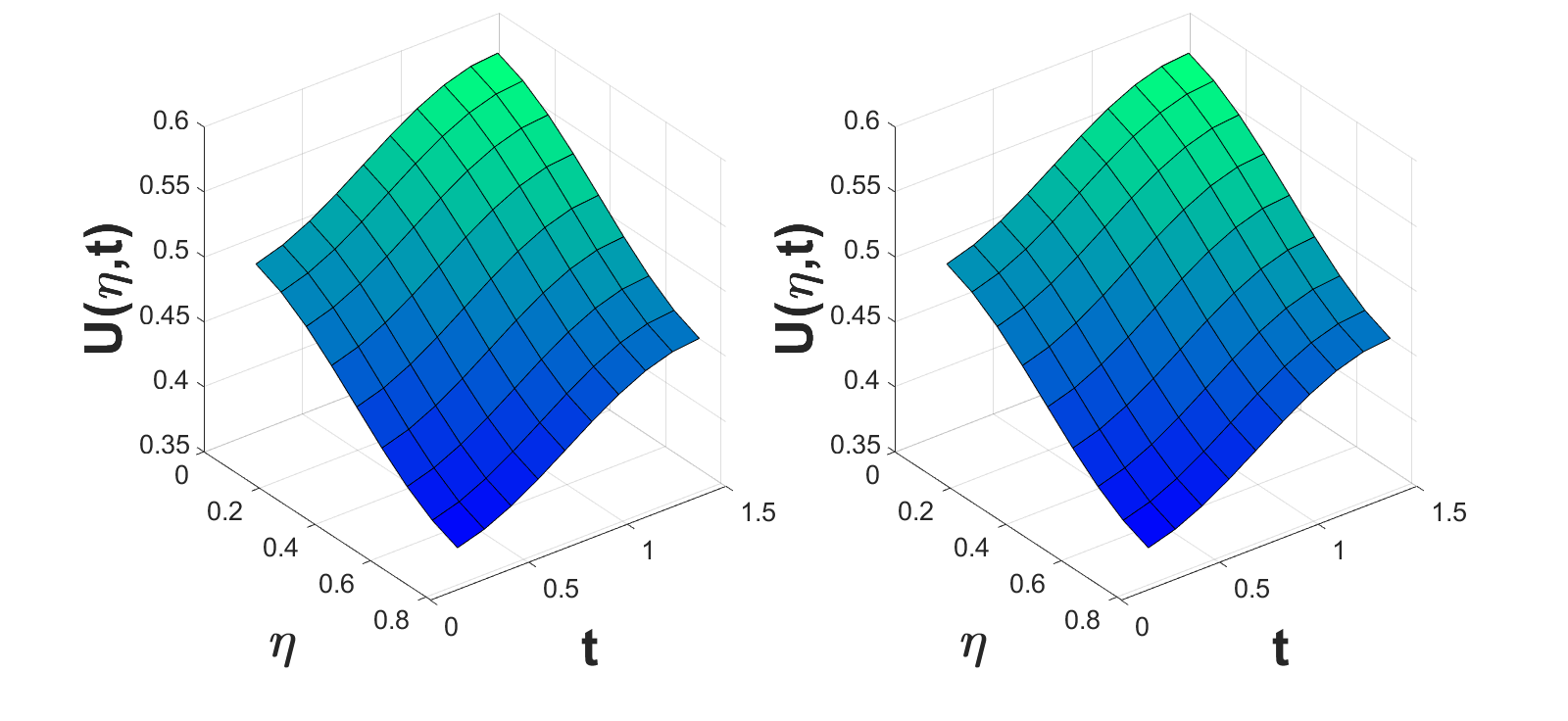} \\
     \includegraphics[width=14cm, height=6.3cm]{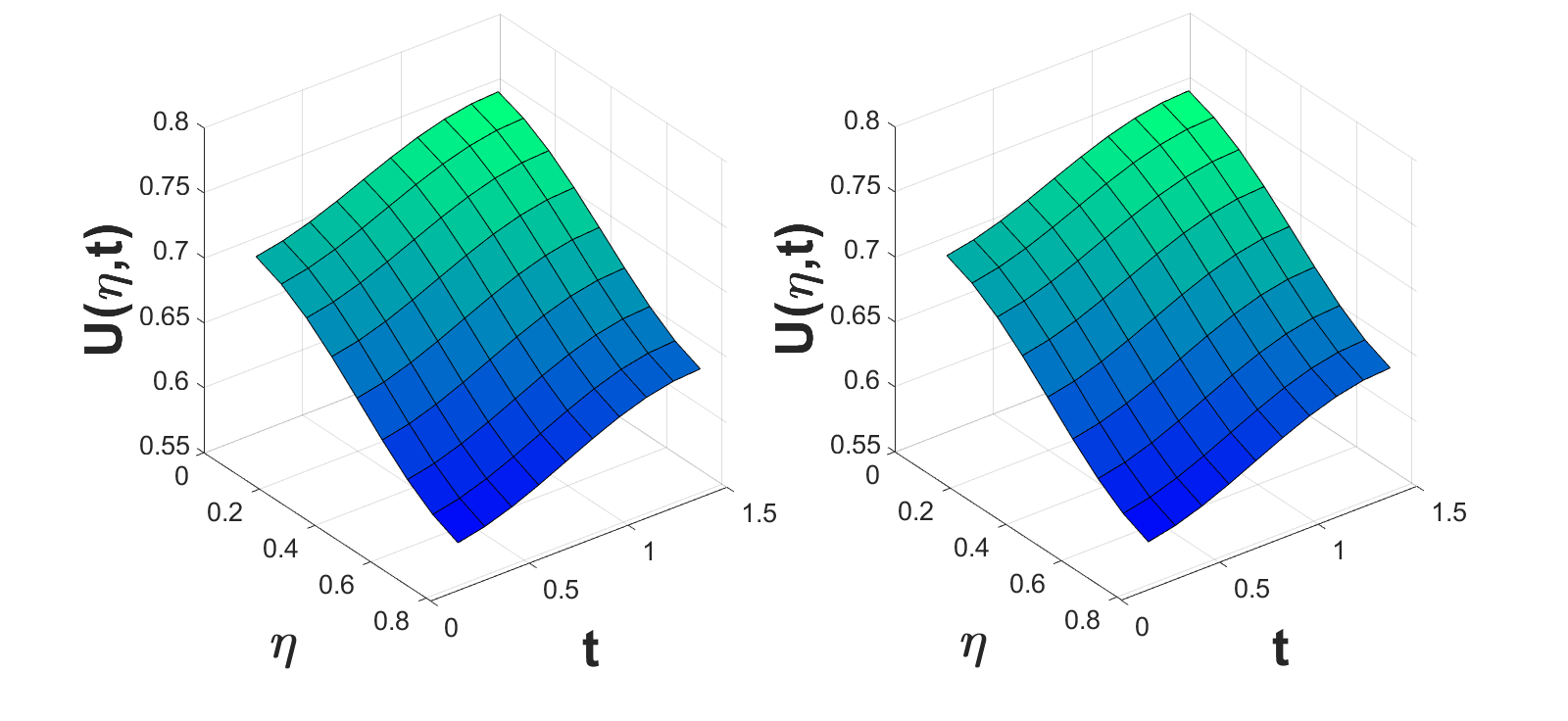}\\ 
      \includegraphics[width=14cm, height=6.3cm]{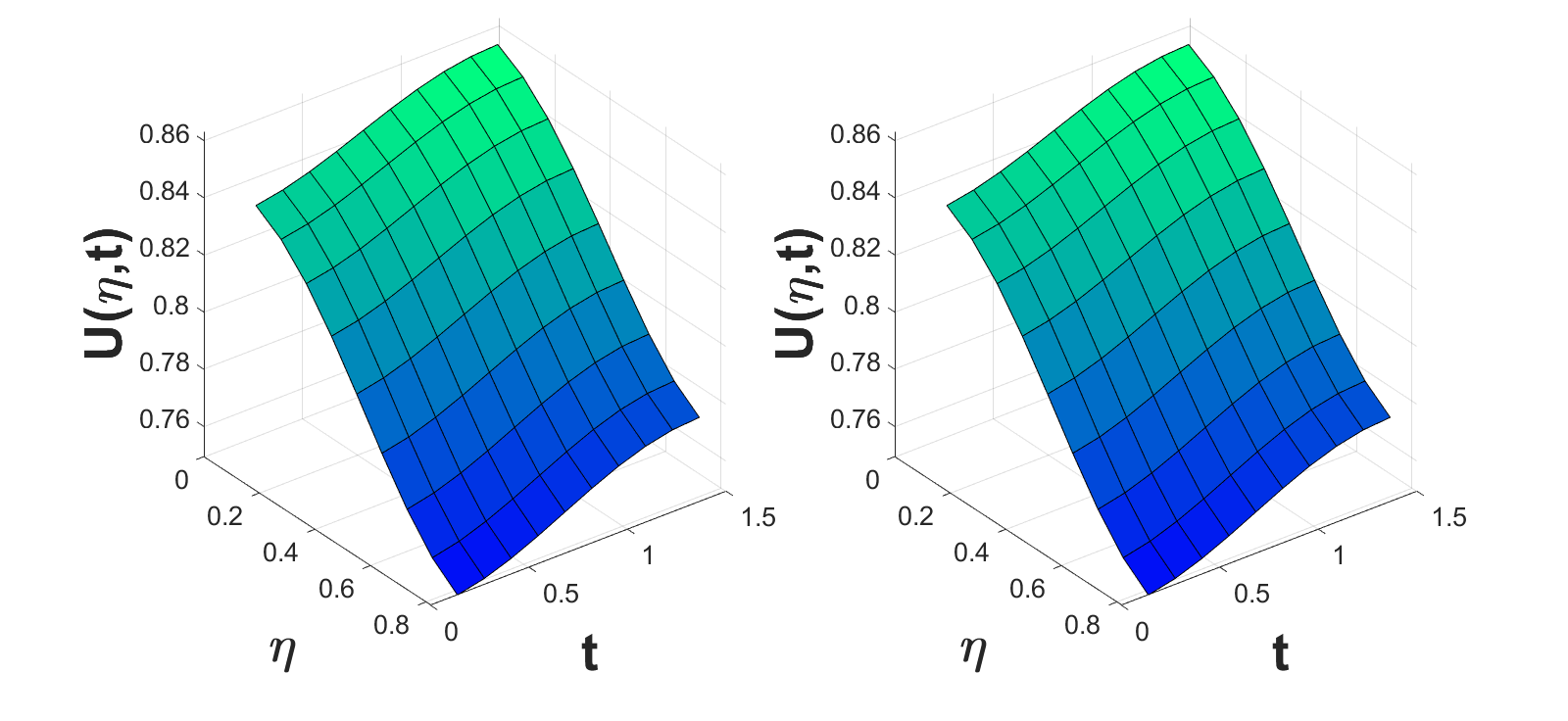}
    \caption{3D graphs for exact(left) vs approximate(right) solution for Example~\ref{eg3} at time $ T=$ 1.5 for $\delta =$ 1, 2, 4, and 8, respectively.}
    \label{fig5_eg3}
  \end{figure}

     \begin{example}\label{eg1}
      \end{example}
  In this example, gBF equation~\ref{Eq.1.4} is approximated for particular small and negative  values $\sigma_1 =0.1$, and $\sigma_2 =-0.0025$ using TS-CPsM. Increasing values of $\delta$ make the equation highly non-linear, which is taken care of by the proposed method with excellent accuracy. At a very small number of grid points, we have obtained highly accurate results at different times in terms of the error infinity norm, which are efficient compared to those existing in the literature~\cite{zhu2010numerical}. The calculated results have been verified for different sets of values of $\delta$ and $t$, and the same is reported in Table~\ref{tab1}. Results show that for a very small number of grid points $N=6$ in both space and time directions, we achieved very good accuracy.
Figure~\ref{fig1_eg1}  shows a comparison of the numerical and exact solutions for $\delta = 1,~2,~4,~\text{and}~ 8$ at various time values. 
A solid line indicates an exact solution, whereas a block indicates a numerical result at a certain time step. Further, 3D graphs for exact vs approximate solution at time $ T=$ 0.3 for $\delta =$ 1, 2, 4, and 8, respectively, are given in Figure~\ref{fig2_eg1}.
 \begin{table}
 \centering
 \caption{\label{tab1}Numerical solution proposed for Example~\ref{eg1} with different time t for various values of $\delta$ and grid points $N$.}
 \begin{tabular}{ ccccc }
       \toprule
                                              
   & t&	N=4  &    N= 6 & \\
   \hline 	
   &	 &  TS-CPsM&         TS-CPsM    & BSQI~\cite{zhu2010numerical} \\
   \midrule
 &0.1&          &  &  \\						
 &$\delta$  = 1& 9.5696e-13&	5.5511e-17&	1.3240e-11\\
 &$\delta$ =2&	1.4242e-12&	4.5519e-15&		2.8470e-10\\
 &$\delta$ =4&	9.9620e-13&	1.1102e-16&	3.9917e-10\\
 &$\delta$=8&	4.7429e-13&	1.1102e-16& -	\\
 &0.2&		&		&		                \\
 &$\delta$=1&	1.0758e-12	&5.5511e-17&	1.7803e-11	\\
 &$\delta$ =2&	1.6033e-12&	2.2204e-16&		3.8795e-10\\
 &$\delta$ =4&	1.1265e-12&	1.1102e-16&	5.4380e-10\\
 &$\delta$=8&	5.4279e-13	&1.1102e-16&	-\\
 						
 &0.3&	&	&		\\
 &$\delta$ =1&	1.1146e-12&	5.5511e-17&	1.9426e-11	\\
 &$\delta$ =2&	1.6627e-12&	1.1102e-16&		4.2465e-10\\
 &$\delta$ =4&	1.1720e-12&	1.1102e-16&		5.9517e-10\\
 &$\delta$=8&	5.6999e-13&	1.1102e-16&		-\\
 						
 &0.4	&	&		&		\\
 &$\delta$  =1&	1.1343e-12&	5.5511e-17&		2.0008e-11\\
 &$\delta$ =2&	1.6921e-12&	1.1102e-16&		4.3759e-10\\
 &$\delta$=4&	1.1960e-12&	1.1102e-16&		6.1323e-10\\
 &$\delta$ =8&	5.8642e-13&	1.1102e-16&	-	\\
 						
 &0.5&	&	&		\\	
 &$\delta$ =1&	1.1478e-12&	5.5511e-17&		2.0216e-11\\
 &$\delta$ =2&	1.7099e-12&	1.1102e-16&	4.0205e-10\\
 &$\delta$ =4&	1.2114e-12&	1.1102e-16&	6.1941e-10\\
 &$\delta$ =8&	5.9841e-13&	1.1102e-16&		-\\
  \bottomrule
 \end{tabular}
 \end{table} 
     \begin{example}\label{eg2}
 \end{example}
 In the second example of the queue the considered Eq.~\ref{Eq.1.4} is approximated  using TS-CPsM for a new set of positive and equal values $\sigma_1 =1$, and $\sigma_2 =1$. As may be seen from Table~\ref{tab2}, with an increase in grid points, the accuracy is getting higher and higher.
Here, we found highly accurate results for 16 grid points in both space and time;
 however, the literature says that a large number of grid points in the temporal direction is required to get good results. A further observation is that there is a slight decay in accuracy as the delta values are increased for a large time, though the results are still comparably better and quite acceptable.
  Graphical outcomes of approximated solutions with the exact solutions for different sets of values of time with varying $\delta$ are shown in Figure~\ref{fig3a_eg2}. A solution that is exact is shown by a solid line, whereas numerical results at discrete times are shown as blocks. Further, 3D graphs for exact vs approximate solution at time $ T=$ 0.6 for $\delta =$ 1, 2, 4, and 8, respectively, are given in Figure~\ref{fig3b_eg2}.
 \begin{sidewaystable}
\centering
\caption{\label{tab2}Numerical solution proposed for Example~\ref{eg2} with different time t for various values of $\delta$ and grid points $N$.}
\begin{tabular}{ ccccccccccc }
      \toprule                                             
  & &	N=4  &                         N=6&	N=8    &              N=10 & N=12& N=14 & N=16 & & \\
  \hline 	
  &	 &  TS-CPsM & 	TS-CPsM   & TS-CPsM &   TS-CPsM  & TS-CPsM&TS-CPsM&TS-CPsM & BSQI~\cite{zhu2010numerical}&SCM~\cite{javidi2006modified}\\
  \midrule
&t=0.2& &	&	&	&	&	& & & \\			
&$\delta$= 1&	1.094e-07&	1.363e-10&	1.665e-13&	2.220e-16&	1.110e-16&	1.110e-16&	1.110e-16&5.558e-7	& 8.749e-12\\
&$\delta$=2&	1.430e-07&	2.267e-10&	3.968e-13&	8.882e-16&	1.110e-16&	1.110e-16&	2.220e-16&2.561e-6	& 5.812e-11\\
&$\delta$ =4&	1.722e-07&	1.145e-09&	6.861e-12&	3.775e-14&	2.220e-16&	1.110e-16&	1.110e-16&1.762e-6	& 1.397e-10\\
&$\delta$=8&	2.574e-06&	8.626e-08&	2.139e-09&	4.481e-11&	8.373e-13&	1.421e-14&	2.220e-16&-	& 2.972e-10\\
&t=0.4&		&	&		&	&	&& 	&&\\
&$\delta$ =1&	1.008e-07&	2.416e-10&	6.424e-13&	1.499e-15&	1.110e-16&	1.110e-16&	1.110e-16&9.055e-7	& 1.667e-11\\
&$\delta$ =2&	4.296e-07&	3.601e-09&	2.289e-11&	1.376e-13&	6.661e-16&	2.220e-16&	2.220e-16&4.243e-6	& 6.388e-11\\
&$\delta$ =4&	3.209e-06&	1.508e-07&	4.579e-09&	1.098e-10&	2.264e-12&	4.152e-14&	6.661e-16&4.174e-7	& 1.365e-10\\
&$\delta$ =8&	5.131e-05&	3.319e-06&	7.566e-08&	1.537e-08&	1.609e-09&	9.949e-11&	4.360e-12&-	& 2.169e-10\\
								
&t=0.6&	&	&	&	&	& &&  &\\			
&$\delta$=1&	3.683e-07&	2.995e-09&	1.816e-11&	9.320e-14&	4.441e-16&	1.110e-16&	1.110e-16&2.188e-6	& 2.274e-11\\
&$\delta$=2&	5.567e-07&	1.280e-08&	4.287e-10&	1.083e-11&	2.0140e-13&	3.220e-15&	2.220e-16&3.569e-6	& 6.190e-11\\
&$\delta$=4&	2.722e-05&	2.085e-06&	9.336e-08&	2.570e-09&	3.230e-11&	3.884e-12&	2.531e-13&2.424e-6	& 9.537e-11\\
&$\delta$ =8&	5.763e-05&	1.164e-05&	3.690e-06&	3.129e-07&	1.654e-08&	5.834e-09&	6.217e-10&-	& 2.722e-11\\
								
&t=0.8&	&	&	&	&	&& 	& & \\		
&$\delta$ =1&	9.945e-07&	1.384e-08&	1.356e-10&	1.047e-12&	6.939e-15&	1.110e-16&	1.110e-16&2.933e-6	& 2.655e-11\\
&$\delta$ =2&	3.882e-06&	2.548e-07&	1.141e-08&	3.577e-10&	8.942e-12&	1.873e-13&	3.331e-15&1.465e-6	& 6.153e-11\\
&$\delta$=4&	6.466e-05&	5.138e-06&	1.557e-07&	3.309e-08&	3.627e-09&	2.205e-10&	7.847e-12&2.358e-6	& 5.230e-11\\
&$\delta$ =8&	8.497e-05&	5.466e-05&	7.401e-06&	1.286e-06&	4.016e-07&	1.778e-08&	8.669e-09	& -&5.573e-12\\ 
								
&t=1.0&	&	&	&& 	&	&	& &\\		
&$\delta$ =1&	1.811e-06&	3.257e-08&	3.246e-10&	2.000e-12&	6.017e-14&	1.499e-15&	1.110e-16&3.015e-6& 	2.813e-11\\
&$\delta$ =2&	1.296e-05&	1.245e-06&	6.765e-08	&2.454e-09	&6.029e-11&	7.717e-13&	3.853e-14	&5.542e-6& 5.320e-11\\
&$\delta$=4&	8.701e-05&	5.130e-06	&1.653e-06	&2.581e-07&	1.563e-08&	6.608e-10&	1.886e-10&	1.444e-6& 2.533e-11\\
&$\delta$ =8&	1.431e-04&	8.879e-05&	9.984e-06&	6.545e-06&	4.035e-07&	2.768e-07&	4.092e-08&-	& 1.111e-12\\ 
 \bottomrule
\end{tabular}
\end{sidewaystable}
     \begin{example}\label{eg3}
     \end{example} 
     In this demonstration, a different combination of the variables $\sigma_1 =1$ and $\sigma_2 = 0$ is taken into account, which results in the transformation of the gBF equation into the more straightforward Burgers equation.  The importance and use of the equation have already been discussed in the introduction part. The error infinity norm for increasing number of grid points is shown in Table~\ref{tab3} for a different set of values of $\delta$ and time $t$. Figure~\ref{fig4_eg3}   shows a comparison of the numerical and exact solutions for $\delta = 1,~2,~\text{and}~4$.
     A solid line indicates an exact solution, whereas a block indicates a numerical result at a certain time step. 
      Further, 3D graphs for exact vs approximate solution at time $ T=1.5$ for $\delta =$ 1, 2, and 4, respectively, is given in Figure~\ref{fig5_eg3}.

\begin{figure}[htbp]
         \centering
        \includegraphics[width=6.48cm, height=6.2cm]{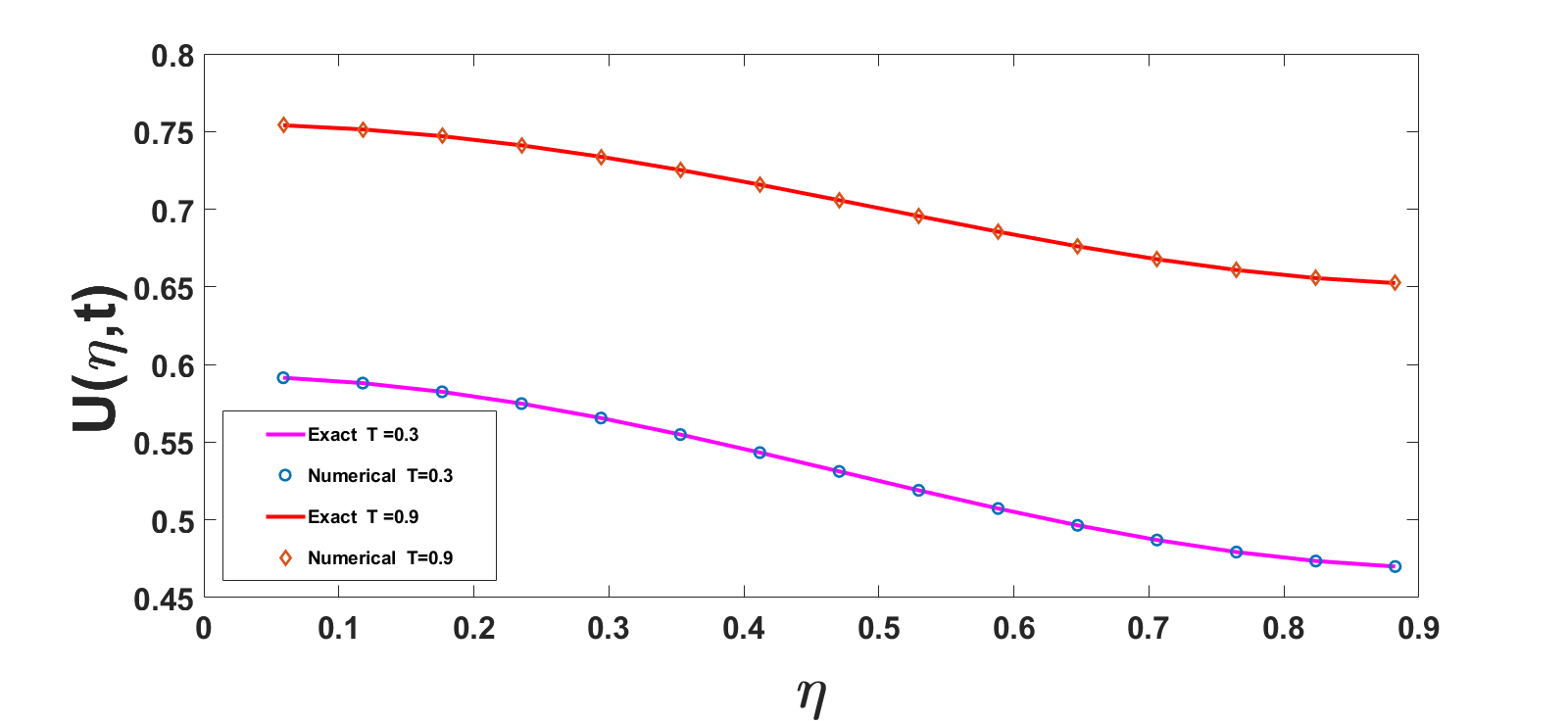} 
       \includegraphics[width=6.48cm, height=6.2cm]{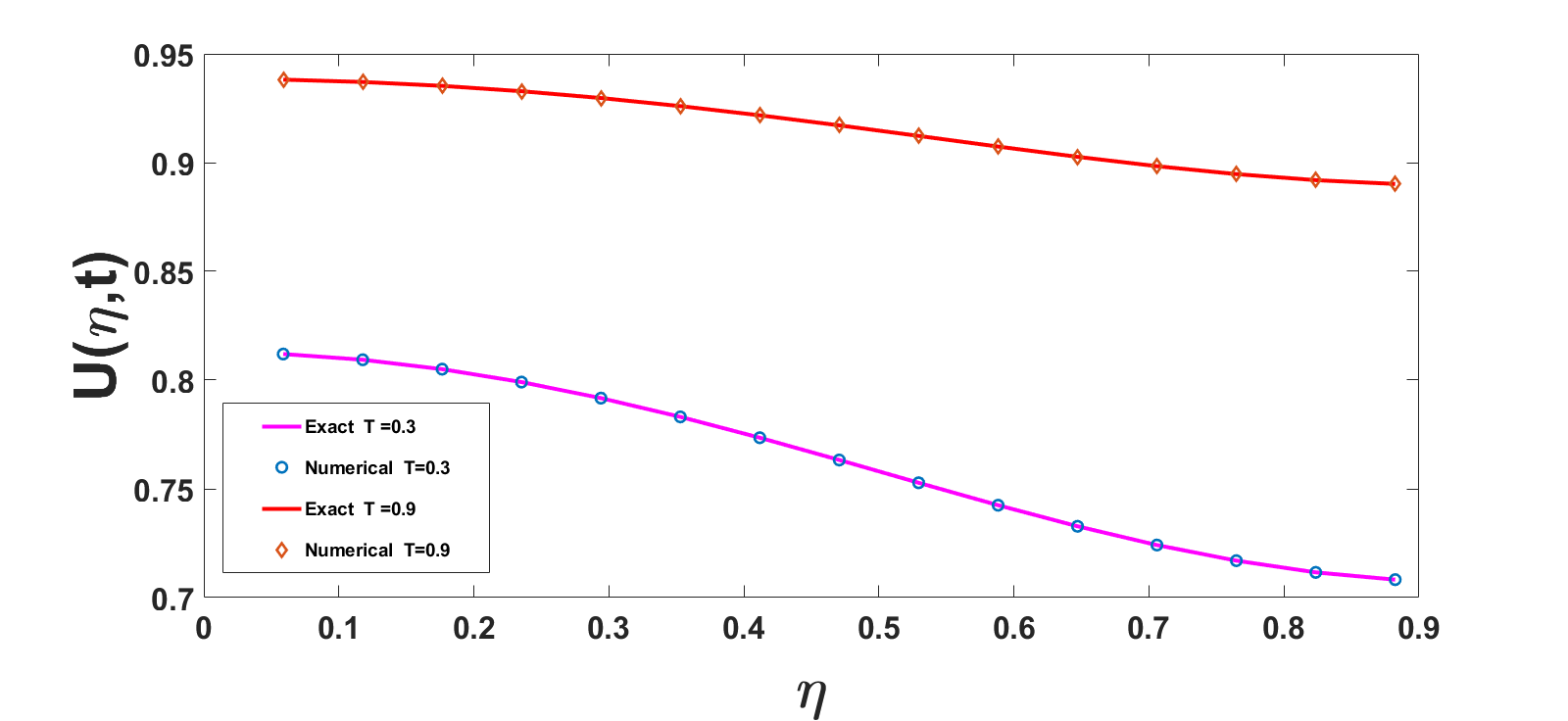}\\
       \includegraphics[width=6.48cm, height=6.2cm]{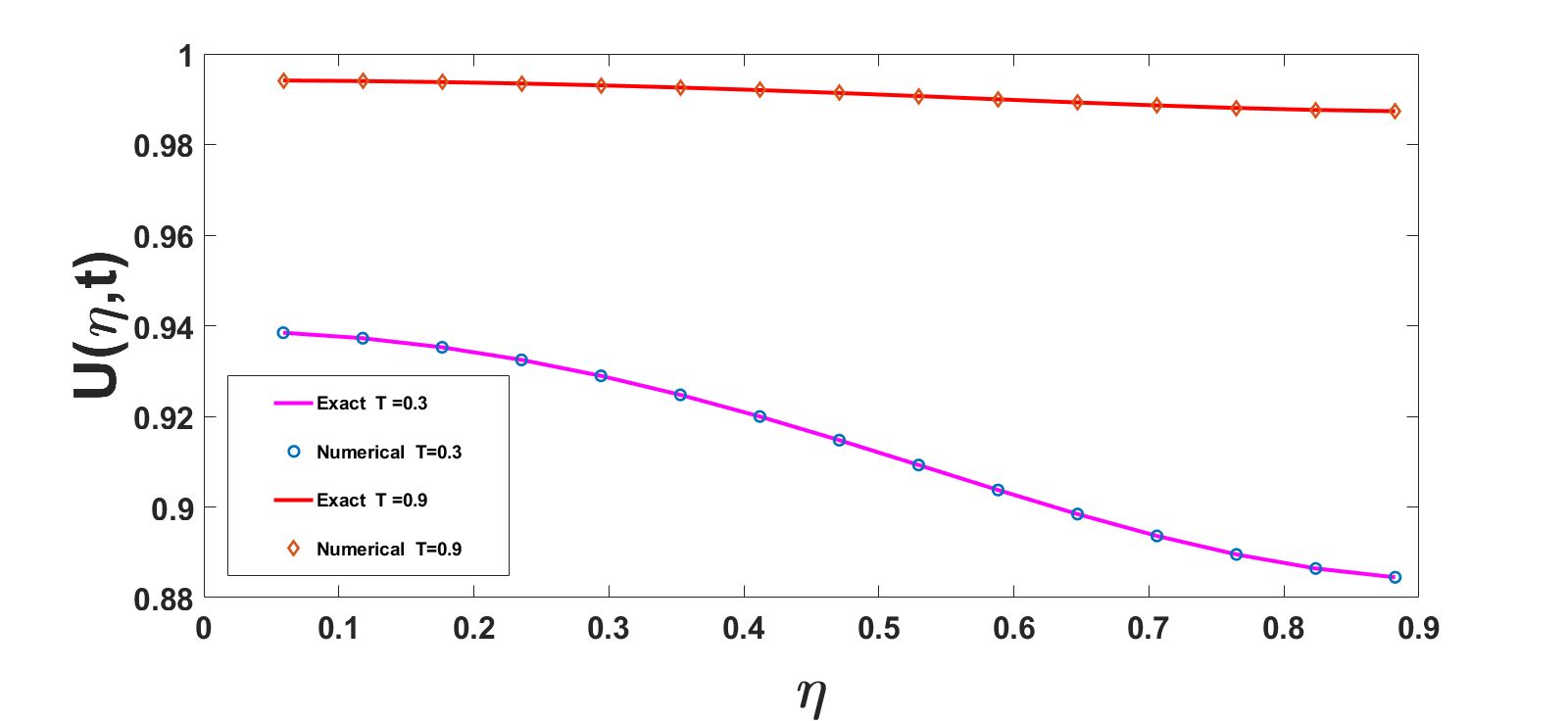}
        \includegraphics[width=6.48cm, height=6.2cm]{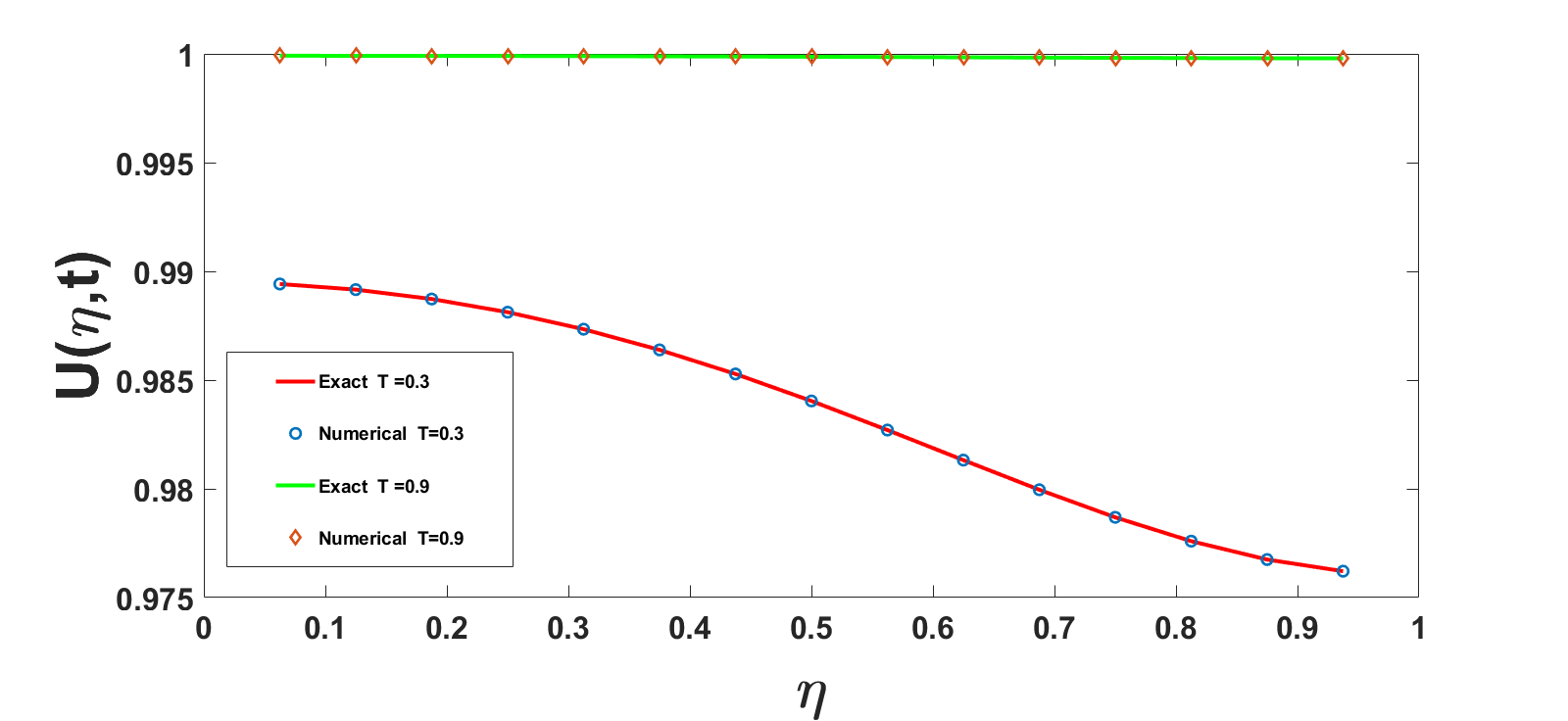}
        \caption{Numerical  vs Exact solution for Example~\ref{eg4} when time $ T=$ 0.3 and 0.9   for $\delta =$ 1, 2, 4, and 8, respectively.}
         \label{fig6_eg4}
        \end{figure}

\begin{table}
\centering
\caption{\label{tab3}Numerical solution proposed for Example~\ref{eg3} with different time t for various values of $\delta$ and grid points $N$.}
\begin{tabular}{ ccccccc }
      \toprule
                                             
  & &	N=4  &                        N=6&	N=8    &               N=10  & \\
  \hline 	
  &	 &  TS-CPsM &     	TS-CPsM   &   	TS-CPsM &          TS-CPsM   &SCM~\cite{javidi2006modified} \\
  \midrule
%

&t=0.001&	&&&&		\\				
&$\delta$=1&	2.7715e-09&	1.1814e-11&	3.5527e-14&	1.1102e-16&	6.5918e-13\\
&$\delta$ =2&	5.8538e-09&	4.3071e-11&	1.9740e-13&	7.7716e-16&		4.0404e-13\\
&$\delta$=4&	5.8482e-09	&6.1648e-11&	4.0001e-13&	1.7764e-15&			1.8770e-13\\

&t=0.5&	&	&	&	&\\		
&$\delta$ =1&	1.1039e-07&	1.3468e-10&	1.5776e-13	&2.2204e-16	&1.7366e-13\\
&$\delta$=2	&1.9814e-07&	3.6723e-10&	7.2031e-13&	1.5543e-15&	4.5848e-13\\
&$\delta$=4&	1.9663e-07&	4.7425e-10&	1.3300e-12&	3.7748e-15&			1.9439e-13\\

&t=1&	&&	& & \\						
&$\delta$=1&	1.1644e-07&	1.4060e-10&	1.6886e-13&	2.2204e-16&		3.5012e-13\\
&$\delta$=2	&1.8895e-07&	3.3341e-10&	6.5703e-13	&1.4433e-15	&	5.0652e-13\\
&$\delta$ =4&	1.8534e-07	&4.4473e-10&	1.2240e-12&	3.6637e-15&			1.9994e-13\\

&t=1.5&	&	&		&	&		\\
&$\delta$=1	&1.1146e-07&	1.3011e-10&	1.5110e-13&	2.2204e-16&		5.1602e-13\\
&$\delta$=2&	1.7174e-07&	2.9088e-10&	5.3968e-13&	1.1102e-15&		5.4759e-13\\
&$\delta$ =4&	1.6998e-07&	4.0300e-10&	1.0719e-12	&3.1086e-15	&		2.0447e-13\\

&t=2	&	&	&		&	&\\
&$\delta$ =1&	1.0571e-07&	1.1338e-10	&2.2149e-13&	6.1062e-16&	7.9885e-13\\
&$\delta$ =2&	1.4803e-07&	2.2619e-10	&3.9735e-13	&7.7716e-15&		5.8118e-13\\
&$\delta$ =4&	1.5113e-07&	3.4518e-10&	8.9762e-13&	2.5535e-15&				2.0795e-13\\

&t=2.5&		&	&	&	&	\\	
&$\delta$=1&	9.9698e-08&	2.3955e-10&	1.2607e-12&	5.5511e-15&		6.6683e-13\\
&$\delta$ =2	&1.1989e-07&	1.7146e-10&	3.2796e-13&	6.6613e-16&			6.0701e-13\\
&$\delta$=4&	1.2930e-07&	2.7381e-10&	6.8279e-13&	1.9984e-15&			2.1038e-13\\
 \bottomrule
\end{tabular}
\end{table}        
  \vspace{1cm} 
 \begin{example}\label{eg4}  
 \end{example}
Consideration of a different set of time values has been found in the literature ~\cite{zhu2010numerical},~\cite{golbabai2009spectral} for the parameters $\sigma_1 =1$ and  $\sigma_2 =1$ and the approximation for the same set is made and compared in this example.
Table~\ref{tab4} displays the error infinity norms for varying $\delta$ and $t$. As can be seen from this table, with an increase in the grid points $N$, the error is decreasing and the results are getting more accurate. A comparison with the existing results shows the efficiency and robustness of the proposed technique. The comparison between the exact and approximated results is depicted in Figure~\ref{fig6_eg4} for varying $\delta$ at various time values $t$. 

\begin{sidewaystable}
\centering
\caption{\label{tab4}Numerical solution proposed for Example~\ref{eg4} with different time t for various values of $\delta$ and grid points $N$.}
\begin{tabular}{ ccccccccccc }
      \toprule
                                             
  & &	N=4  &    &	&N=8    &  & &N=16 && \\
  \hline 	
  &	 &  TS-CPsM &     BSQI~\cite{zhu2010numerical} & SDD~\cite{golbabai2009spectral} &  TS-CPsM &  BSQI~\cite{zhu2010numerical}& SDD~\cite{golbabai2009spectral}& TS-CPsM& BSQI~\cite{zhu2010numerical} &SDD~\cite{golbabai2009spectral}\\
  \midrule 
&t=0.3&	&	&	&	&	&&&&	\\
&$\delta$ =1&	1.096e-07&	1.350e-6& 1.476e-5&	1.503e-13&	5.773e-6& 4.165e-11&	1.110e-16&	2.218e-6&1.184e-10\\
&$\delta$ =2&1.665e-07	& -&-	&3.332e-12	&	-	&-&2.220e-16 &-&-	\\
&$\delta$ =4&	4.247e-07&	1.230e-5&9.159e-6	&1.765e-10&	5.750e-5&1.285e-10&	1.110e-16&	3.274e-5&2.119e-10\\
&$\delta$ =8&	2.133e-05&	1.638e-6&3.956e-6	&4.988e-08&	3.919e-6&8.851e-12&	7.383e-14&	4.151e-7&2.119e-10\\
&t=0.9&	&	&	&	&	&&&&\\	
&$\delta$ =1&	1.400e-06&	9.770e-6&4.110e-5&	2.494e-10&	4.384e-5&9.668e-11&	1.110e-16&	2.480e-5&5.109e-10\\
&$\delta$ =2&7.276e-06	& -&-	&3.124e-08	&	-	&-&8.549e-15 &-&-	\\
&$\delta$ =4&	7.917e-05	&3.017e-6&1.189e-5&	5.254e-07&	1.342e-5&3.891e-11&	2.268e-11&	7.292e-6&7.681e-10\\
&$\delta$ =8&	1.058e-04&	1.116e-7&4.334e-6&	4.914e-06&	4.923e-7&6.158e-12&	2.730e-08&	2.843e-7&7.016e-10\\
 \bottomrule
\end{tabular}
\end{sidewaystable}

\newpage
\section{Conclusion}\label{sec5}
Using the time-space Chebyshev pseudospectral method (TS-CPsM), numerical solutions for the gBF equation with various values of advection and reaction coefficients are reported with tremendous accuracy in this study. It has been observed that for a very small number of grid points $N$ and for large values of non-linear factor $\delta$, the approximated results for the considered equation using the suggested technique have been proven to be highly competent with exact solutions and comparably much better with those existing in the literature, including ~\cite{kaya2004numerical},~\cite{zhu2010numerical},~\cite{ismail2004adomian},~\cite{javidi2006modified},~\cite{ismail2004restrictive}, and others referenced.
The proposed method's stability analysis has been performed. Four different cases for the considered problem have been discussed in order to demonstrate the effectiveness and reliability of the method, and all of them give excellent and encouraging numerical results.
\bmhead{Acknowledgments}
The Ministry of Human Resource Development of India provided financial assistance for this study, which the first author gratefully acknowledges.
\bmhead{Conflict of interest}
The authors declare no conflict of interests.\\

\bibliography{sn-article}


\begin{thebibliography}{51}
\ifx \bisbn   \undefined \def \bisbn  #1{ISBN #1}\fi
\ifx \binits  \undefined \def \binits#1{#1}\fi
\ifx \bauthor  \undefined \def \bauthor#1{#1}\fi
\ifx \batitle  \undefined \def \batitle#1{#1}\fi
\ifx \bjtitle  \undefined \def \bjtitle#1{#1}\fi
\ifx \bvolume  \undefined \def \bvolume#1{\textbf{#1}}\fi
\ifx \byear  \undefined \def \byear#1{#1}\fi
\ifx \bissue  \undefined \def \bissue#1{#1}\fi
\ifx \bfpage  \undefined \def \bfpage#1{#1}\fi
\ifx \blpage  \undefined \def \blpage #1{#1}\fi
\ifx \burl  \undefined \def \burl#1{\textsf{#1}}\fi
\ifx \doiurl  \undefined \def \doiurl#1{\url{https://doi.org/#1}}\fi
\ifx \betal  \undefined \def \betal{\textit{et al.}}\fi
\ifx \binstitute  \undefined \def \binstitute#1{#1}\fi
\ifx \binstitutionaled  \undefined \def \binstitutionaled#1{#1}\fi
\ifx \bctitle  \undefined \def \bctitle#1{#1}\fi
\ifx \beditor  \undefined \def \beditor#1{#1}\fi
\ifx \bpublisher  \undefined \def \bpublisher#1{#1}\fi
\ifx \bbtitle  \undefined \def \bbtitle#1{#1}\fi
\ifx \bedition  \undefined \def \bedition#1{#1}\fi
\ifx \bseriesno  \undefined \def \bseriesno#1{#1}\fi
\ifx \blocation  \undefined \def \blocation#1{#1}\fi
\ifx \bsertitle  \undefined \def \bsertitle#1{#1}\fi
\ifx \bsnm \undefined \def \bsnm#1{#1}\fi
\ifx \bsuffix \undefined \def \bsuffix#1{#1}\fi
\ifx \bparticle \undefined \def \bparticle#1{#1}\fi
\ifx \barticle \undefined \def \barticle#1{#1}\fi
\bibcommenthead
\ifx \bconfdate \undefined \def \bconfdate #1{#1}\fi
\ifx \botherref \undefined \def \botherref #1{#1}\fi
\ifx \url \undefined \def \url#1{\textsf{#1}}\fi
\ifx \bchapter \undefined \def \bchapter#1{#1}\fi
\ifx \bbook \undefined \def \bbook#1{#1}\fi
\ifx \bcomment \undefined \def \bcomment#1{#1}\fi
\ifx \oauthor \undefined \def \oauthor#1{#1}\fi
\ifx \citeauthoryear \undefined \def \citeauthoryear#1{#1}\fi
\ifx \endbibitem  \undefined \def \endbibitem {}\fi
\ifx \bconflocation  \undefined \def \bconflocation#1{#1}\fi
\ifx \arxivurl  \undefined \def \arxivurl#1{\textsf{#1}}\fi
\csname PreBibitemsHook\endcsname

\bibitem[\protect\citeauthoryear{Bateman}{1915}]{bateman1915some}
\begin{barticle}
\bauthor{\bsnm{Bateman}, \binits{H.}}:
\batitle{Some recent researches on the motion of fluids}.
\bjtitle{Monthly Weather Review}
\bvolume{43}(\bissue{4}),
\bfpage{163}--\blpage{170}
(\byear{1915})
\end{barticle}
\endbibitem

\bibitem[\protect\citeauthoryear{Burgers}{1948}]{BURGERS1948171}
\begin{botherref}
\oauthor{\bsnm{Burgers}, \binits{J.M.}}:
A mathematical model illustrating the theory of turbulence.
Advances in Applied Mechanics,
vol. 1,
pp. 171--199.
Elsevier
(1948)
\end{botherref}
\endbibitem

\bibitem[\protect\citeauthoryear{Greenshields
  et~al.}{1935}]{greenshields1935study}
\begin{bchapter}
\bauthor{\bsnm{Greenshields}, \binits{B.D.}},
\bauthor{\bsnm{Bibbins}, \binits{J.}},
\bauthor{\bsnm{Channing}, \binits{W.}},
\bauthor{\bsnm{Miller}, \binits{H.}}:
\bctitle{A study of traffic capacity}.
In: \bbtitle{Highway Research Board Proceedings},
vol. \bseriesno{14},
pp. \bfpage{448}--\blpage{477}
(\byear{1935}).
\bcomment{Washington, DC}
\end{bchapter}
\endbibitem

\bibitem[\protect\citeauthoryear{Murray}{1973}]{murray1973burgers}
\begin{barticle}
\bauthor{\bsnm{Murray}, \binits{J.}}:
\batitle{On burgers' model equations for turbulence}.
\bjtitle{Journal of Fluid Mechanics}
\bvolume{59}(\bissue{2}),
\bfpage{263}--\blpage{279}
(\byear{1973})
\end{barticle}
\endbibitem

\bibitem[\protect\citeauthoryear{Hirsh}{1975}]{hirsh1975higher}
\begin{barticle}
\bauthor{\bsnm{Hirsh}, \binits{R.S.}}:
\batitle{Higher order accurate difference solutions of fluid mechanics problems
  by a compact differencing technique}.
\bjtitle{Journal of computational physics}
\bvolume{19}(\bissue{1}),
\bfpage{90}--\blpage{109}
(\byear{1975})
\end{barticle}
\endbibitem

\bibitem[\protect\citeauthoryear{Yepez}{2002}]{yepez2002efficient}
\begin{botherref}
\oauthor{\bsnm{Yepez}, \binits{J.}}:
An efficient quantum algorithm for the one-dimensional burgers equation.
arXiv preprint quant-ph/0210092
(2002)
\end{botherref}
\endbibitem

\bibitem[\protect\citeauthoryear{Fisher}{1937}]{fisher1937wave}
\begin{barticle}
\bauthor{\bsnm{Fisher}, \binits{R.A.}}:
\batitle{The wave of advance of advantageous genes}.
\bjtitle{Annals of eugenics}
\bvolume{7}(\bissue{4}),
\bfpage{355}--\blpage{369}
(\byear{1937})
\end{barticle}
\endbibitem

\bibitem[\protect\citeauthoryear{Bramson}{1978}]{bramson1978maximal}
\begin{barticle}
\bauthor{\bsnm{Bramson}, \binits{M.D.}}:
\batitle{Maximal displacement of branching brownian motion}.
\bjtitle{Communications on Pure and Applied Mathematics}
\bvolume{31}(\bissue{5}),
\bfpage{531}--\blpage{581}
(\byear{1978})
\end{barticle}
\endbibitem

\bibitem[\protect\citeauthoryear{Britton et~al.}{1986}]{britton1986reaction}
\begin{bbook}
\bauthor{\bsnm{Britton}, \binits{N.F.}}, \betal:
\bbtitle{Reaction-diffusion Equations and Their Applications to Biology.}
\bpublisher{Academic Press}, \blocation{???}
(\byear{1986})
\end{bbook}
\endbibitem

\bibitem[\protect\citeauthoryear{Frank-Kamenetskii}{2015}]{frank2015diffusion}
\begin{bbook}
\bauthor{\bsnm{Frank-Kamenetskii}, \binits{D.A.}}:
\bbtitle{Diffusion and Heat Exchange in Chemical Kinetics}
vol. \bseriesno{2171}.
\bpublisher{Princeton University Press}, \blocation{???}
(\byear{2015})
\end{bbook}
\endbibitem

\bibitem[\protect\citeauthoryear{Canosa}{1969}]{canosa1969diffusion}
\begin{barticle}
\bauthor{\bsnm{Canosa}, \binits{J.}}:
\batitle{Diffusion in nonlinear multiplicative media}.
\bjtitle{Journal of Mathematical Physics}
\bvolume{10}(\bissue{10}),
\bfpage{1862}--\blpage{1868}
(\byear{1969})
\end{barticle}
\endbibitem

\bibitem[\protect\citeauthoryear{Ablowitz et~al.}{1987}]{ablowitz1987topics}
\begin{bbook}
\bauthor{\bsnm{Ablowitz}, \binits{M.J.}},
\bauthor{\bsnm{Fuchssteiner}, \binits{B.}},
\bauthor{\bsnm{Kruskal}, \binits{M.}}:
\bbtitle{Topics in Soliton Theory and Exactly Solvable Nonlinear Equations:
  Proceedings of the Conference on Nonlinear Evolution Equations, Solitons and
  the Inverse Scattering Transform}.
\bpublisher{World Scientific}, \blocation{???}
(\byear{1987})
\end{bbook}
\endbibitem

\bibitem[\protect\citeauthoryear{Wang and Lu}{1990}]{wang1990exact}
\begin{barticle}
\bauthor{\bsnm{Wang}, \binits{X.}},
\bauthor{\bsnm{Lu}, \binits{Y.}}:
\batitle{Exact solutions of the extended burgers-fisher equation}.
\bjtitle{Chinese Physics Letters}
\bvolume{7}(\bissue{4}),
\bfpage{145}--\blpage{147}
(\byear{1990})
\end{barticle}
\endbibitem

\bibitem[\protect\citeauthoryear{Chandrasekaran and
  Ramasami}{1996}]{chandrasekaran1996painleve}
\begin{barticle}
\bauthor{\bsnm{Chandrasekaran}, \binits{P.}},
\bauthor{\bsnm{Ramasami}, \binits{E.}}:
\batitle{Painleve analysis of a class of nonlinear diffusion equations}.
\bjtitle{International Journal of Stochastic Analysis}
\bvolume{9}(\bissue{1}),
\bfpage{77}--\blpage{86}
(\byear{1996})
\end{barticle}
\endbibitem

\bibitem[\protect\citeauthoryear{Chen and Zhang}{2004}]{chen2004new}
\begin{barticle}
\bauthor{\bsnm{Chen}, \binits{H.}},
\bauthor{\bsnm{Zhang}, \binits{H.}}:
\batitle{New multiple soliton solutions to the general burgers--fisher equation
  and the kuramoto--sivashinsky equation}.
\bjtitle{Chaos, Solitons \& Fractals}
\bvolume{19}(\bissue{1}),
\bfpage{71}--\blpage{76}
(\byear{2004})
\end{barticle}
\endbibitem

\bibitem[\protect\citeauthoryear{Lu et~al.}{2007}]{lu2007some}
\begin{barticle}
\bauthor{\bsnm{Lu}, \binits{J.}},
\bauthor{\bsnm{Yu-Cui}, \binits{G.}},
\bauthor{\bsnm{Shu-Jiang}, \binits{X.}}:
\batitle{Some new exact solutions to the burgers--fisher equation and
  generalized burgers--fisher equation}.
\bjtitle{Chinese Physics}
\bvolume{16}(\bissue{9}),
\bfpage{2514}
(\byear{2007})
\end{barticle}
\endbibitem

\bibitem[\protect\citeauthoryear{Kaya and El-Sayed}{2004}]{kaya2004numerical}
\begin{barticle}
\bauthor{\bsnm{Kaya}, \binits{D.}},
\bauthor{\bsnm{El-Sayed}, \binits{S.M.}}:
\batitle{A numerical simulation and explicit solutions of the generalized
  burgers--fisher equation}.
\bjtitle{Applied Mathematics and computation}
\bvolume{152}(\bissue{2}),
\bfpage{403}--\blpage{413}
(\byear{2004})
\end{barticle}
\endbibitem

\bibitem[\protect\citeauthoryear{Zhu and Kang}{2010}]{zhu2010numerical}
\begin{barticle}
\bauthor{\bsnm{Zhu}, \binits{C.-G.}},
\bauthor{\bsnm{Kang}, \binits{W.-S.}}:
\batitle{Numerical solution of burgers--fisher equation by cubic b-spline
  quasi-interpolation}.
\bjtitle{Applied Mathematics and Computation}
\bvolume{216}(\bissue{9}),
\bfpage{2679}--\blpage{2686}
(\byear{2010})
\end{barticle}
\endbibitem

\bibitem[\protect\citeauthoryear{Mittal and
  Tripathi}{2015}]{mittal2015numerical}
\begin{barticle}
\bauthor{\bsnm{Mittal}, \binits{R.}},
\bauthor{\bsnm{Tripathi}, \binits{A.}}:
\batitle{Numerical solutions of generalized burgers--fisher and generalized
  burgers--huxley equations using collocation of cubic b-splines}.
\bjtitle{International Journal of Computer Mathematics}
\bvolume{92}(\bissue{5}),
\bfpage{1053}--\blpage{1077}
(\byear{2015})
\end{barticle}
\endbibitem

\bibitem[\protect\citeauthoryear{Singh et~al.}{2020}]{singh2020fourth}
\begin{barticle}
\bauthor{\bsnm{Singh}, \binits{A.}},
\bauthor{\bsnm{Dahiya}, \binits{S.}},
\bauthor{\bsnm{Singh}, \binits{S.}}:
\batitle{A fourth-order b-spline collocation method for nonlinear
  burgers--fisher equation}.
\bjtitle{Mathematical Sciences}
\bvolume{14},
\bfpage{75}--\blpage{85}
(\byear{2020})
\end{barticle}
\endbibitem

\bibitem[\protect\citeauthoryear{Rashidi et~al.}{2009}]{rashidi2009explicit}
\begin{barticle}
\bauthor{\bsnm{Rashidi}, \binits{M.}},
\bauthor{\bsnm{Ganji}, \binits{D.}},
\bauthor{\bsnm{Dinarvand}, \binits{S.}}:
\batitle{Explicit analytical solutions of the generalized burger and
  burger--fisher equations by homotopy perturbation method}.
\bjtitle{Numerical Methods for Partial Differential Equations: An International
  Journal}
\bvolume{25}(\bissue{2}),
\bfpage{409}--\blpage{417}
(\byear{2009})
\end{barticle}
\endbibitem

\bibitem[\protect\citeauthoryear{Ismail et~al.}{2004}]{ismail2004adomian}
\begin{barticle}
\bauthor{\bsnm{Ismail}, \binits{H.N.}},
\bauthor{\bsnm{Raslan}, \binits{K.}},
\bauthor{\bsnm{Abd~Rabboh}, \binits{A.A.}}:
\batitle{Adomian decomposition method for burger's--huxley and burger's--fisher
  equations}.
\bjtitle{Applied mathematics and computation}
\bvolume{159}(\bissue{1}),
\bfpage{291}--\blpage{301}
(\byear{2004})
\end{barticle}
\endbibitem

\bibitem[\protect\citeauthoryear{Khattak}{2009}]{khattak2009computational}
\begin{barticle}
\bauthor{\bsnm{Khattak}, \binits{A.J.}}:
\batitle{A computational meshless method for the generalized burger's--huxley
  equation}.
\bjtitle{Applied Mathematical Modelling}
\bvolume{33}(\bissue{9}),
\bfpage{3718}--\blpage{3729}
(\byear{2009})
\end{barticle}
\endbibitem

\bibitem[\protect\citeauthoryear{Xu and Xian}{2010}]{xu2010application}
\begin{barticle}
\bauthor{\bsnm{Xu}, \binits{Z.-h.}},
\bauthor{\bsnm{Xian}, \binits{D.-q.}}:
\batitle{Application of exp-function method to generalized burgers-fisher
  equation}.
\bjtitle{Acta Mathematicae Applicatae Sinica, English Series}
\bvolume{26}(\bissue{4}),
\bfpage{669}--\blpage{676}
(\byear{2010})
\end{barticle}
\endbibitem

\bibitem[\protect\citeauthoryear{Zhang and Yan}{2010}]{zhang2010lattice}
\begin{barticle}
\bauthor{\bsnm{Zhang}, \binits{J.}},
\bauthor{\bsnm{Yan}, \binits{G.}}:
\batitle{A lattice boltzmann model for the burgers--fisher equation}.
\bjtitle{Chaos: An Interdisciplinary Journal of Nonlinear Science}
\bvolume{20}(\bissue{2}),
\bfpage{023129}
(\byear{2010})
\end{barticle}
\endbibitem

\bibitem[\protect\citeauthoryear{Sari et~al.}{2010}]{sari2010compact}
\begin{barticle}
\bauthor{\bsnm{Sari}, \binits{M.}},
\bauthor{\bsnm{G{\"u}rarslan}, \binits{G.}},
\bauthor{\bsnm{Da{\u{g}}}, \binits{{\. I}.}}:
\batitle{A compact finite difference method for the solution of the generalized
  burgers--fisher equation}.
\bjtitle{Numerical Methods for Partial Differential Equations: An International
  Journal}
\bvolume{26}(\bissue{1}),
\bfpage{125}--\blpage{134}
(\byear{2010})
\end{barticle}
\endbibitem

\bibitem[\protect\citeauthoryear{Sari}{2011}]{sari2011differential}
\begin{barticle}
\bauthor{\bsnm{Sari}, \binits{M.}}:
\batitle{Differential quadrature solutions of the generalized burgers--fisher
  equation with a strong stability preserving high-order time integration}.
\bjtitle{Mathematical and Computational Applications}
\bvolume{16}(\bissue{2}),
\bfpage{477}--\blpage{486}
(\byear{2011})
\end{barticle}
\endbibitem

\bibitem[\protect\citeauthoryear{Zhang et~al.}{2012}]{zhang2012local}
\begin{barticle}
\bauthor{\bsnm{Zhang}, \binits{R.}},
\bauthor{\bsnm{Yu}, \binits{X.}},
\bauthor{\bsnm{Zhao}, \binits{G.}}:
\batitle{The local discontinuous galerkin method for burger's--huxley and
  burger's--fisher equations}.
\bjtitle{Applied Mathematics and Computation}
\bvolume{218}(\bissue{17}),
\bfpage{8773}--\blpage{8778}
(\byear{2012})
\end{barticle}
\endbibitem

\bibitem[\protect\citeauthoryear{Nawaz et~al.}{2013}]{nawaz2013application}
\begin{botherref}
\oauthor{\bsnm{Nawaz}, \binits{R.}},
\oauthor{\bsnm{Ullah}, \binits{H.}},
\oauthor{\bsnm{Islam}, \binits{S.}},
\oauthor{\bsnm{Idrees}, \binits{M.}}:
Application of optimal homotopy asymptotic method to burger equations.
Journal of Applied Mathematics
\textbf{2013}
(2013)
\end{botherref}
\endbibitem

\bibitem[\protect\citeauthoryear{Yadav and Jiwari}{2017}]{yadav2017finite}
\begin{barticle}
\bauthor{\bsnm{Yadav}, \binits{O.P.}},
\bauthor{\bsnm{Jiwari}, \binits{R.}}:
\batitle{Finite element analysis and approximation of burgers'-fisher
  equation}.
\bjtitle{Numerical Methods for Partial Differential Equations}
\bvolume{33}(\bissue{5}),
\bfpage{1652}--\blpage{1677}
(\byear{2017})
\end{barticle}
\endbibitem

\bibitem[\protect\citeauthoryear{Malik et~al.}{2015}]{malik2015numerical}
\begin{barticle}
\bauthor{\bsnm{Malik}, \binits{S.A.}},
\bauthor{\bsnm{Qureshi}, \binits{I.M.}},
\bauthor{\bsnm{Amir}, \binits{M.}},
\bauthor{\bsnm{Malik}, \binits{A.N.}},
\bauthor{\bsnm{Haq}, \binits{I.}}:
\batitle{Numerical solution to generalized burgers'-fisher equation using
  exp-function method hybridized with heuristic computation}.
\bjtitle{PLoS One}
\bvolume{10}(\bissue{3}),
\bfpage{0121728}
(\byear{2015})
\end{barticle}
\endbibitem

\bibitem[\protect\citeauthoryear{Nadeem et~al.}{2019}]{nadeem2019modified}
\begin{barticle}
\bauthor{\bsnm{Nadeem}, \binits{M.}},
\bauthor{\bsnm{Li}, \binits{F.}},
\bauthor{\bsnm{Ahmad}, \binits{H.}}:
\batitle{Modified laplace variational iteration method for solving fourth-order
  parabolic partial differential equation with variable coefficients}.
\bjtitle{Computers \& Mathematics with Applications}
\bvolume{78}(\bissue{6}),
\bfpage{2052}--\blpage{2062}
(\byear{2019})
\end{barticle}
\endbibitem

\bibitem[\protect\citeauthoryear{Loyinmi and Akinfe}{2020}]{loyinmi2020exact}
\begin{barticle}
\bauthor{\bsnm{Loyinmi}, \binits{A.C.}},
\bauthor{\bsnm{Akinfe}, \binits{T.K.}}:
\batitle{Exact solutions to the family of fisher's reaction-diffusion equation
  using elzaki homotopy transformation perturbation method}.
\bjtitle{Engineering Reports}
\bvolume{2}(\bissue{2}),
\bfpage{12084}
(\byear{2020})
\end{barticle}
\endbibitem

\bibitem[\protect\citeauthoryear{Akinfe and Loyinmi}{2021}]{akinfe2021solitary}
\begin{barticle}
\bauthor{\bsnm{Akinfe}, \binits{T.K.}},
\bauthor{\bsnm{Loyinmi}, \binits{A.C.}}:
\batitle{A solitary wave solution to the generalized burgers-fisher's equation
  using an improved differential transform method: A hybrid scheme approach}.
\bjtitle{Heliyon}
\bvolume{7}(\bissue{5}),
\bfpage{07001}
(\byear{2021})
\end{barticle}
\endbibitem

\bibitem[\protect\citeauthoryear{Mohanty and Sharma}{2020}]{mohanty2020high}
\begin{barticle}
\bauthor{\bsnm{Mohanty}, \binits{R.K.}},
\bauthor{\bsnm{Sharma}, \binits{S.}}:
\batitle{A high-resolution method based on off-step non-polynomial spline
  approximations for the solution of burgers-fisher and coupled nonlinear
  burgers' equations}.
\bjtitle{Engineering Computations}
\bvolume{37}(\bissue{8}),
\bfpage{2785}--\blpage{2818}
(\byear{2020})
\end{barticle}
\endbibitem

\bibitem[\protect\citeauthoryear{Kumar~Verma and
  Kayenat}{2019}]{kumar2019stability}
\begin{barticle}
\bauthor{\bsnm{Kumar~Verma}, \binits{A.}},
\bauthor{\bsnm{Kayenat}, \binits{S.}}:
\batitle{On the stability of micken's type nsfd schemes for generalized burgers
  fisher equation}.
\bjtitle{Journal of Difference Equations and Applications}
\bvolume{25}(\bissue{12}),
\bfpage{1706}--\blpage{1737}
(\byear{2019})
\end{barticle}
\endbibitem

\bibitem[\protect\citeauthoryear{Li}{2019}]{li2019geometric}
\begin{barticle}
\bauthor{\bsnm{Li}, \binits{J.}}:
\batitle{Geometric properties and exact travelling wave solutions for the
  generalized burger-fisher equation and the sharma-tasso-olver equation}.
\bjtitle{Journal of Nonlinear Modeling and analysis}
\bvolume{1}(\bissue{1}),
\bfpage{1}--\blpage{10}
(\byear{2019})
\end{barticle}
\endbibitem

\bibitem[\protect\citeauthoryear{Onyejekwe
  et~al.}{2020}]{onyejekwe2020numerical}
\begin{barticle}
\bauthor{\bsnm{Onyejekwe}, \binits{O.O.}},
\bauthor{\bsnm{Minale}, \binits{B.}},
\bauthor{\bsnm{Habtamu}, \binits{F.}},
\bauthor{\bsnm{Amha}, \binits{T.}},
\bauthor{\bsnm{Tamiru}, \binits{G.}},
\bauthor{\bsnm{Mengistu}, \binits{B.}},
\bauthor{\bsnm{Demiss}, \binits{Y.}},
\bauthor{\bsnm{Alemseged}, \binits{N.}}, \betal:
\batitle{Numerical discrete-domain integral formulations for generalized
  burger-fisher equation}.
\bjtitle{Applied Mathematics}
\bvolume{11}(\bissue{03}),
\bfpage{137}
(\byear{2020})
\end{barticle}
\endbibitem

\bibitem[\protect\citeauthoryear{Zhang et~al.}{2021}]{zhang2021global}
\begin{barticle}
\bauthor{\bsnm{Zhang}, \binits{H.}},
\bauthor{\bsnm{Xia}, \binits{Y.}},
\bauthor{\bsnm{N'gbo}, \binits{P.-R.}}:
\batitle{Global existence and uniqueness of a periodic wave solution of the
  generalized burgers--fisher equation}.
\bjtitle{Applied Mathematics Letters}
\bvolume{121},
\bfpage{107353}
(\byear{2021})
\end{barticle}
\endbibitem

\bibitem[\protect\citeauthoryear{Mittal and Balyan}{2021}]{mittal2021numerical}
\begin{barticle}
\bauthor{\bsnm{Mittal}, \binits{A.K.}},
\bauthor{\bsnm{Balyan}, \binits{L.K.}}:
\batitle{Numerical solutions of time and space fractional coupled burgers
  equations using time--space chebyshev pseudospectral method}.
\bjtitle{Mathematical Methods in the Applied Sciences}
\bvolume{44}(\bissue{4}),
\bfpage{3127}--\blpage{3137}
(\byear{2021})
\end{barticle}
\endbibitem

\bibitem[\protect\citeauthoryear{Mittal}{2020}]{mittal2020spectrally}
\begin{barticle}
\bauthor{\bsnm{Mittal}, \binits{A.}}:
\batitle{Spectrally accurate approximate solutions and convergence analysis of
  fractional burgers' equation}.
\bjtitle{Arabian Journal of Mathematics}
\bvolume{9}(\bissue{3}),
\bfpage{633}--\blpage{644}
(\byear{2020})
\end{barticle}
\endbibitem

\bibitem[\protect\citeauthoryear{Balyan et~al.}{2020}]{balyan2020stability}
\begin{barticle}
\bauthor{\bsnm{Balyan}, \binits{L.}},
\bauthor{\bsnm{Mittal}, \binits{A.}},
\bauthor{\bsnm{Kumar}, \binits{M.}},
\bauthor{\bsnm{Choube}, \binits{M.}}:
\batitle{Stability analysis and highly accurate numerical approximation of
  fisher's equations using pseudospectral method}.
\bjtitle{Mathematics and Computers in Simulation}
\bvolume{177},
\bfpage{86}--\blpage{104}
(\byear{2020})
\end{barticle}
\endbibitem

\bibitem[\protect\citeauthoryear{Saini et~al.}{2022a}]{saini2022comparative}
\begin{botherref}
\oauthor{\bsnm{Saini}, \binits{P.}},
\oauthor{\bsnm{Balyan}, \binits{L.}},
\oauthor{\bsnm{Kumar}, \binits{A.}},
\oauthor{\bsnm{Singh}, \binits{G.}}:
Comparative analysis of post-processing on spectral collocation methods for
  non-smooth functions.
Signal, Image and Video Processing,
1--9
(2022)
\end{botherref}
\endbibitem

\bibitem[\protect\citeauthoryear{Saini et~al.}{2022b}]{saini2022modification}
\begin{barticle}
\bauthor{\bsnm{Saini}, \binits{P.}},
\bauthor{\bsnm{Balyan}, \binits{L.}},
\bauthor{\bsnm{Kumar}, \binits{A.}},
\bauthor{\bsnm{Singh}, \binits{G.}}:
\batitle{Modification of chebyshev pseudospectral method to minimize the gibbs
  oscillatory behaviour in resynthesizing process}.
\bjtitle{Circuits, Systems, and Signal Processing}
\bvolume{41}(\bissue{11}),
\bfpage{6238}--\blpage{6265}
(\byear{2022})
\end{barticle}
\endbibitem

\bibitem[\protect\citeauthoryear{Hesthaven
  et~al.}{2007}]{hesthaven2007spectral}
\begin{bbook}
\bauthor{\bsnm{Hesthaven}, \binits{J.S.}},
\bauthor{\bsnm{Gottlieb}, \binits{S.}},
\bauthor{\bsnm{Gottlieb}, \binits{D.}}:
\bbtitle{Spectral Methods for Time-dependent Problems}
vol. \bseriesno{21}.
\bpublisher{Cambridge University Press}, \blocation{???}
(\byear{2007})
\end{bbook}
\endbibitem

\bibitem[\protect\citeauthoryear{Boyd}{2001}]{boyd2001chebyshev}
\begin{bbook}
\bauthor{\bsnm{Boyd}, \binits{J.P.}}:
\bbtitle{Chebyshev and Fourier Spectral Methods}.
\bpublisher{Courier Corporation}, \blocation{???}
(\byear{2001})
\end{bbook}
\endbibitem

\bibitem[\protect\citeauthoryear{Trefethen}{2000}]{trefethen2000spectral}
\begin{barticle}
\bauthor{\bsnm{Trefethen}, \binits{L.N.}}:
\batitle{Spectral methods in matlab, volume 10 of software, environments, and
  tools}.
\bjtitle{Society for Industrial and Applied Mathematics (SIAM), Philadelphia,
  PA}
\bvolume{24},
\bfpage{57}
(\byear{2000})
\end{barticle}
\endbibitem

\bibitem[\protect\citeauthoryear{Canuto et~al.}{2012}]{canuto2012others}
\begin{botherref}
\oauthor{\bsnm{Canuto}, \binits{C.}},
\oauthor{\bsnm{Hussaini}, \binits{M.}},
\oauthor{\bsnm{Quarteroni}, \binits{A.}},
\oauthor{\bsnm{Thomas~Jr}, \binits{A.}}:
others. Spectral methods in fluid dynamics.
Springer
(2012)
\end{botherref}
\endbibitem

\bibitem[\protect\citeauthoryear{Javidi}{2006}]{javidi2006modified}
\begin{bchapter}
\bauthor{\bsnm{Javidi}, \binits{M.}}:
\bctitle{Modified pseudospectral method for generalized burger's-fisher
  equation}.
In: \bbtitle{International Mathematical Forum},
vol. \bseriesno{1},
pp. \bfpage{1555}--\blpage{1564}
(\byear{2006}).
\bcomment{Citeseer}
\end{bchapter}
\endbibitem

\bibitem[\protect\citeauthoryear{Golbabai and
  Javidi}{2009}]{golbabai2009spectral}
\begin{barticle}
\bauthor{\bsnm{Golbabai}, \binits{A.}},
\bauthor{\bsnm{Javidi}, \binits{M.}}:
\batitle{A spectral domain decomposition approach for the generalized
  burger's--fisher equation}.
\bjtitle{Chaos, Solitons \& Fractals}
\bvolume{39}(\bissue{1}),
\bfpage{385}--\blpage{392}
(\byear{2009})
\end{barticle}
\endbibitem

\bibitem[\protect\citeauthoryear{Ismail and
  Abd~Rabboh}{2004}]{ismail2004restrictive}
\begin{barticle}
\bauthor{\bsnm{Ismail}, \binits{H.N.}},
\bauthor{\bsnm{Abd~Rabboh}, \binits{A.A.}}:
\batitle{A restrictive pad{\'e} approximation for the solution of the
  generalized fisher and burger--fisher equations}.
\bjtitle{Applied Mathematics and Computation}
\bvolume{154}(\bissue{1}),
\bfpage{203}--\blpage{210}
(\byear{2004})
\end{barticle}
\endbibitem

\end{thebibliography}

\end{document}